\documentclass[11pt]{article}
\usepackage{geometry}                % See geometry.pdf to learn the layout options. There are lots.
\usepackage{graphicx}
\usepackage{authblk}
\usepackage{amssymb}
\usepackage{epstopdf}
\usepackage[sans]{dsfont}
\usepackage[T1]{fontenc}
\usepackage[english]{babel}
\usepackage{latexsym}
\usepackage{mathrsfs}
\usepackage{amscd}
\usepackage{graphicx}
\usepackage{color}
\usepackage{float}
\usepackage{enumitem}
\usepackage{amsmath}
\usepackage{amsfonts}
\numberwithin{equation}{section}
\usepackage{enumerate}
\usepackage{amsthm}
\usepackage[numbers]{natbib}
\usepackage[bookmarks=true,colorlinks=true,linkcolor={blue},urlcolor={blue}, citecolor={blue},pdfstartview={XYZ null null 1.22}]{hyperref}%
\graphicspath{{pics/}}

\usepackage{geometry}                
\usepackage{epstopdf}

\newtheorem{thm}{Theorem}[section]
\newtheorem{remark}{Remark}

\newtheorem{prop}[thm]{Proposition}

\definecolor{Red}{rgb}{1, 0, 0} % Red of svgnames
\newcommand{\be}{\begin{equation}}
\newcommand{\ee}{\end{equation}}

\newcommand{\epsi}{\varepsilon}
\newcommand{\ba}{\begin{array}}
\newcommand{\ea}{\end{array}}
\newcommand{\pM}{\left(\begin{array}}
\newcommand{\Mp}{\end{array}\right)}

\newcommand{\R}{\mathbb{R}}

\newcommand{\La}{\mathcal{L}}
\newcommand{\I}{\mathcal{I}}

\begin{document}
%%%%%%%%%%%%%%%%%%%%%%%%%%%%%%%%%%%%%%%%%%%%%%%
\title{Kinetic description of optimal control problems and applications to opinion consensus}
\author[1]{Giacomo Albi\thanks{giacomo.albi@unife.it}}
\author[2]{Michael Herty\thanks{herty@igpm.rwth-aachen.de}}
\author[1]{Lorenzo Pareschi\thanks{lorenzo.pareschi@unife.it}}
\affil[1]{\em University of Ferrara, Department of Mathematics and Computer Science, Via Machiavelli 35, I-44121 Ferrara, Italy}
\affil[2]{\em RWTH Aachen University, Department of Mathematics, Templergraben 55, D-52065 Aachen, Germany}

\date{\today}
\maketitle

\begin{abstract}
In this paper an optimal control problem for a large system of interacting agents is considered using a kinetic perspective. As a prototype model  we analyze a microscopic model
of opinion formation under constraints.  For this problem a Boltzmann--type equation based  on a model predictive control formulation is introduced and discussed. In particular, the receding horizon  strategy permits to embed the minimization of suitable cost functional into  binary particle interactions. 
The corresponding Fokker-Planck asymptotic limit is also derived and explicit expressions of stationary solutions are given. Several numerical results showing the robustness of the present approach are finally reported. 
\end{abstract}

%\tableofcontents
\section{Introduction}
The development of mathematical models describing the collective behavior of systems of interacting agents originated a large literature in the recent years with applications to several fields, like biology, engineering, economy and sociology (see \cite{BAT:13,BS12,CordierPareschiToscani2005aa,CKFL:07,CuSm:07,DCBC:06,GGS:82,HK02,LuMa:99,Sznajd00,Tos06,HertyRinghofer2011ad,HertyRinghofer2011ac,ArmbrusterRinghofer2005aa,HertyPareschi2010aa} and the references therein). Most of these models are at the level of the microscopic dynamic described by a system of ordinary differential equations. Only recently some of these models have been related to partial differential equations through the corresponding kinetic and hydrodynamic description \cite{AlbiPareschi2013aa,ArmbrusterRinghofer2005aa,BS:09,DLMP:13,DegondMotsch2008ab,FoHaTo:11,GGL:12,HaTa:08,HertyRinghofer2011ad,HertyPareschi2010aa,MaPa:12,Tos06}.
We refer to the recent surveys in \cite{MoTa:13,NPT:10,ViZa:12} and to the book \cite{PT:13} for an introduction to the subject.

In this paper we consider problems where the collective behavior corresponds to the process of alignment, like in the opinion formation dynamic. Different to the classical approach where individuals are assumed to freely interact with each other, here we are particularly interested in such problems in a constrained setting. We consider feedback type controls for the resulting process and present a kinetic modeling 
including those controls. This can be used to study the exterior influence of the system dynamics to enforce emergence of non spontaneous desired asymptotic states. Classical examples are given by persuading voters to vote for a specific candidate or by influencing buyers towards a given good or asset \cite{Ben:05,DMPW:09,LuMa:99,MaPa:12}. %We also refer to \cite{ColomboPogodaev2012aa} for an example in biology. 
In our model, the external intervention is introduced as an additional control  subject to certain bounds, representing the limitations, in terms of economic resources, media availability, etc., of the opinion maker.

Control mechanisms of self-organized systems have been studied for macroscopic models in \cite{ColomboPogodaev2012aa,ColomboPogodaev2013aa} and for kinetic and hydrodynamic models in \cite{AlbiPareschi2013aa,DMPW:09,HertyRinghofer2011ac}. However, in the above references, the control is modeled as a leader dynamics. Therefore, it is given a priori and represented by a supplementary differential model. Also, in \cite{HertyRinghofer2011ac} the control is modeled a posteriori on the level of the kinetic equation mimicking a classical LQR control approach. Recently, the control of emergent behaviors in multiagent systems has been studied in \cite{CFPT13,FS13} where  the authors develop the idea of {sparse optimization} (for sparse control it is meant that the policy maker intervenes the minimal amount of times on the minimal amount of individual agents) at the microscopic and kinetic level. We refer also to \cite{BFY:13} for results concerning the control of mean-field type systems.  
Contrary to all those approaches we derive a controller using the model predicitive control framework on the microscopic level
and study the related kinetic description for large number of agents. In this way we do not need to prescribe control dynamics a priori or a posteriori but these are obtained automatically based only on the underlying microscopic interactions and a suitable cost functional.

The starting point of our modeling is a general framework which embed several type of collective alignment models. We consider the evolution of $N$ agents where each agent has an opinion $w_i=w_i(t)  \in \I$, $\I=[-1,1]$, $i=1,\ldots,N$ and this opinion can change over time according to
\begin{equation}\label{eq:pbm}
\begin{aligned}
&\dot{w}_i=\frac{1}{N}\sum_{j=1}^{N}P(w_i,w_j)(w_j-w_i)+u,\qquad\qquad\qquad w_i(0) = w_{0i}, 
\end{aligned}
\end{equation}
where the control $u=u(t)$ is given by the minimization of the cost functional over a certain time horizon $T$
\begin{equation}\label{eq:pbc}
\begin{aligned}
&u=\textrm{argmin}\int_0^T\frac1{N}\sum_{j=1}^N\left(\frac{1}{2}(w_j-w_d)^2+ \frac{\nu}{2} u^2\right)ds,\qquad u(t)\in[u_L,u_R].
\end{aligned}
\end{equation}
In the formulation \eqref{eq:pbc} the value $w_d$ is the desired state and $\nu>0$ is a regularization parameter. We chose a least--square
type cost functional for simplicity but other costs can be treated similarly.  We additionally prescribe box 
constraints on the pointwise values of $u(t)$ given by the constants  $u_L$ and $u_R>u_L.$ The bound constraints on $u(t)$ are required in order
 to preserve the bounds for $w_i.$ 
 The dynamic in \eqref{eq:pbm} describes an average process of alignment between the opinions $w_i$  of the $N$ agents. Typically, the function  $P(w,v)$ is such that $0\leq P(w,v) \leq 1$ and represents a measure of the inclination of the agents to change their opinion. Usually such function $P$ follows the assumption that extreme  opinions are more difficult to be influenced by others \cite{GGL:12,Sznajd00,Tos06}. 
Problem \eqref{eq:pbm}-\eqref{eq:pbc} may be reformulated as Mayer's problem and solved by Pontryagin's maximum principle \cite{Sontag1998aa} or dynamic programming. The main drawback of this approach relies on the fact that the equation for the adjoint variable has to be solved backwards in time  over
the full time interval $[0,T]$. In particular,  for large values of $N$ the computational effort  therefore renders the problem unsolvable. Also, an approach $u=\mathcal{P}(x)$ where $\mathcal{P}$ fulfills a Riccati differential equation cannot be pursued here due to the  large dimension of $\mathcal{P} \in \mathbb{R}^{N \times N}$ and a possible general nonlinearity in $P$. This approach is known as LQR controller in the 
engineering literature \cite{krstic1995}.  A standard methodology, when dealing with such complex system, is based on model predictive control where instead of solving the above control problem over the whole time horizon, the system is approximated by an iterative solution over a sequence of finite time steps \cite{CaBo:04,MayneRawlingsRao2000aa,MichalskaMayne1993aa}.  

In order to decrease the complexity of the model when the number of agents is large, a possible approach is to rely on a kinetic description of the process. Along this line of thought, in this work we introduce a Boltzmann model describing the microscopic model in the model predictive control formulation. Moreover, a Fokker-Planck model is derived in the so called quasi-invariant opinion limit. % (corresponding to the mean-field limit of the system). 
The kinetic models presented in this paper share some common features with the Boltzmann model introduced in \cite{Tos06} in the unconstrained case and with the mean-field constrained models in \cite{CFPT13,FS13}. Here, however, a remarkable difference with respect to \cite{CFPT13,FS13} is that, thanks to the receding horizon strategy, the minimization of the cost functional is embedded into the particle interactions. Similarly to \cite{Tos06}, this permits to compute explicitly the stationary solutions of the resulting constrained dynamic. 

The rest of the manuscript is organized as follows. In the next Section we introduce the model predictive control formulation of system \eqref{eq:pbm}-\eqref{eq:pbc}. In Section 3 a binary dynamic corresponding to the constrained system is introduced and a the main properties of the resulting Boltzmann-type kinetic equation are discussed. In particular, estimates for the convergence of the solution towards the desired state are given. Section 4 is devoted to the derivation of the Fokker-Planck model and the computation of explicit stationary solutions for the resulting kinetic equation. Some modeling variants are discussed in Section 5. Finally, in Section 6 several numerical results are reported showing the robustness of the present approach. Some conclusions and future research directions are made at the end of the manuscript.

\section{Model predictive control}
In this section we adapt the idea of the moving horizon controller (or instantaneous control) to derive a computable control $u$ at any time $t.$ Compared with the solution to \eqref{eq:pbm}-\eqref{eq:pbc}
this control will in general only be suboptimal.  Rigorous results on  the properties of $u$ for quadratic cost functional and linear and nonlinear dynamics  are available, for example, in \cite{CaBo:04,MayneMichalska1990aa}. The model predicitive control framework applied here is also called receding horizon strategy or instantaneous control in the engineering literature. 

\subsection{A receding horizon strategy}
We consider a receding horizon strategy with horizon of a single time interval. Hence, instead of solving \eqref{eq:pbm}-\eqref{eq:pbc} on $[0,T]$, we proceeds as follows:
\begin{itemize}[label=$\bullet$]
\item Split the time interval $[0,T]$ in $M$ time intervals of length $\Delta t$ and let  $t^n=\Delta t \; n$.
\item We assume that the control is piecewise constant on time intervals of length $\Delta t>0$,
$$u(t)=\sum\limits_{n=0}^{M-1} u^n \chi_{[t^{n},t^{n+1}]}(t).$$
\item Determine the  value of the control $u^n \in \mathbb{R}$  by solving for a state $\bar{w}_i$ the (reduced) optimization problem                                    
\begin{equation}\label{eq:pbic}
\begin{aligned}
&\dot{w}_i=\frac{1}{N}\sum_{j=1}^{N}P(w_i,w_j)(w_j-w_i)+u,\qquad\qquad\qquad w_i(t^n) = \bar{w}_i, \\
&u^n=\textrm{argmin}_{ u \in \mathbb{R} } \int_{t^n}^{t^{n+1}} \frac1{N}\sum_{j=1}^N\left(\frac{1}{2}(w_j-w_d)^2+ \frac{\nu}{2} u^2\right)ds,\qquad u \in[u_L,u_R].
\end{aligned}
\end{equation}
\item Having the control  $u^n$ on the interval $[t^n, t^{n+1}]$, evolve $w_i$ according to the dynamics 
\begin{align}\label{eq:controldynamics}
\dot{w}_i=\frac{1}{N}\sum_{j=1}^{N}P(w_i,w_j)(w_j-w_i)+u^n
\end{align}
to obtain  the new state $\bar{w}_i = w_i(t^{n+1}).$ 
\item We again solve \eqref{eq:pbic} to obtain $u^{n+1}$ with the modified initial data. 
\item Repeat this procedure until we reach $n \Delta t = T.$ 
 \end{itemize}    
The advantage compared with the problem  \eqref{eq:pbm}-\eqref{eq:pbc} is the reduced complexity of \eqref{eq:pbic} being an optimization problem in a single real--valued 
variable $u^n.$  Furthermore, for the quadratic cost and a suitable discretization of \eqref{eq:controldynamics} the solution to \eqref{eq:pbic} allows an explicit representation of $u^n$ in 
terms of $\bar{w_i}$ and $w_i(t^{n+1})$ provided $u_L=-\infty$ and $u_R=\infty.$  As shown in section \ref{sec2} this allows to reformulate  the previous algorithm as a feedback controlled system which 
in discretized form reads 
\begin{subequations}
\begin{align}
& w^{n+1}_i=w^n_i+\frac{\Delta t}{N}\sum_{j=1}^{N}P^n_{ij}(w^n_j-w^n_i)+\Delta t u^n,\qquad w^0_i= w_{0i},\label{eq:Dfwd}\\
& u^n = -\frac{\Delta t}{\nu N} \sum_{j=1}^N(w_j^{n+1}-w_d).\label{eq:Dcontrol}
\end{align}
\end{subequations}
%Note that similar to the receding horizon controller for linear quadratic problems the action of the control, $u^n=u(t^n)$, 
%\begin{equation}\label{control law}
%u(t^n) = u^n = -\frac{(\Delta t)^2}{\nu N} \sum_{j=1}^N(w_j^{n+1}-%w_d)
%\end{equation}
% is on the scale $\Delta t$. This is in accordance with the continuous case which would 
% require the discretization of a Riccati equation on $[t^n, t^{n+1}]$ with terminal state zero. 

\begin{remark}
Later on, bounds on the control $u$ as in \eqref{eq:pbic} are required in order to guarantee  that  opinions $w_i \in \mathcal{I}$ for all times. Instead
of considering the constrained problem \eqref{eq:pbic} we will present a condition on $\nu$ ensuring this property in the case of a binary interaction model in Proposition \ref{pr:1} below. 
This allows to treat \eqref{eq:pbic} as an unconstrained problem and  does not require to a priori prescribe bounds $u_L$ and $u_R.$ 
Also note that in general the expression of the control $u$ in terms of $w_i^{n+1}$ and $w_d$ as in equation \eqref{eq:Dcontrol} 
would be much more involved if the bound constraints $u_L,u_R$ are present.
\end{remark}

\subsection{Derivation of the feedback controller} \label{sec2}

%In this section we derive the receding horizon control and establish equation \eqref{eq:Dfwd}-\eqref{eq:Dcontrol}. 
We assume for now that $u_L=-\infty$ and $u_R=+\infty$ and  assume sufficient regularity conditions such that any 
minimizer $u\equiv u^{n} \in \mathbb{R}$ to problem \eqref{eq:pbic}  fulfills the necessary first order optimality conditions. 
We further assume that those conditions are also sufficient for optimality and refer to \cite{Sontag1998aa} for more details.

The optimality conditions on $[t^{n}, t^{n+1}]$ and for $\bar{w_i}=w(t^{n})_i$ are given by the set of the  following equations. 
\begin{align*}
\Delta t \; \nu u &= - \frac{1}N \sum\limits_{i=1}^{N} \int_{t^{n}}^{t^{n+1}} \lambda_i dt, \\
\dot{w}_i & =\frac{1}{N}\sum_{j=1}^{N}P(w_i,w_j)(w_j-w_i)+u, \; w_i(t^{n}) = \bar{w_i}, \\
\dot{\lambda}_i &= - (w_i - w_d)  - \frac{1}N \sum\limits_{j=1}^{N}  R_{ij}  , \; \lambda_i(t^{n+1})=0, \\
R_{ij} &= \lambda_i \partial_{w_i}  \left\{   P(w_i,w_j)(w_j-w_i) \right\} + \lambda_j \partial_{w_j} \{ P(w_j,w_i)(w_i-w_j) \}. 
\end{align*} 

The function $\lambda^{n}_i=\lambda_i(t)$ is the (Lagrange) multiplier.  If we
 discretize the adjoint equation (backwards in time) by the implicit Euler scheme we obtain due to the boundary conditions
$$ \lambda^{n}_i = - \Delta t \;  (w^{n+1}_i - w_d) $$
Further, we may solve for $u$ after discretizing the  integral as $\int_{t^{n}}^{t^{n+1}} f(t) dt = \Delta t \; f^{n}$  to obtain 
$$ u = - \frac{ \Delta t \;}{N \nu } \sum\limits_{i=1}^{N}   (w^{n+1}_i - w_d)$$

Applying an explicit Euler discretization to the dynamics for $w_i$ on the time interval $[t^{n}, t^{n+1}]$ and substituting the control we
obtained, we observe  that the feedback control $u$ is given by \eqref{eq:Dcontrol} and hence the  final equation is given by \eqref{eq:Dfwd}. 

The previous derivation is obtained by first computing the continuous optimality system and then applying a suitable discretization.  However, applying first an explicit
Euler discretization and then computing the discrete optimality system leads to the same result.  Indeed, consider the discretization of 
 \eqref{eq:pbm}--\eqref{eq:pbc} in the interval $[t^{n}, t^{n+1}]$ for constant control $u$ and with  $P^n_{ij}=P(w^n_i,w^n_j)$:
\begin{equation}\label{eq:pb ic disc}
\begin{aligned}
&w^{n+1}_i=w^n_i+\frac{\Delta t}{N}\sum_{j=1}^{N}P^n_{ij}(w^n_j-w_i^{n})+\Delta tu,\quad\qquad\qquad w^n_i = \bar{w}_i, \\
&u =\textrm{argmin} ~\frac{\Delta t}{N} \sum_{j=1}^N\left(\frac{1}{2}(w^n_j-w_d)^2+ \frac{\nu}{2} (u^n)^2\right),
\end{aligned}
\end{equation}
The discrete Lagrangian is given by 
\begin{equation}\label{eq:DLag}
\begin{aligned}
\La(w,\lambda,u)=&\Delta t\left(\frac{1}{N}\sum_{k=1}^N(w_k^{n+1}-w_d)^2+\frac{\nu}{2} u^2\right)+\frac{1}{N}\sum_{k=1}^N\lambda_k^0(w_k^n- \bar{w}_k)
\\&+\frac{1}{N}\sum_{k=1}^{N}\lambda^{n+1}_k\left(w^n_k-w^{n+1}_k+\frac{\Delta t}{N}\sum_{j=1}^{N}P^n_{kj}(w^n_j-w^n_k)+\Delta tu \right)
\end{aligned}
\end{equation}
A minimizer to equation \eqref{eq:pb ic disc} fulfills under suitable regularity assumptions the equations \eqref{eq:pb ic disc}, \eqref{adjoint} and \eqref{optimality}.
\begin{align}\label{adjoint} 
\lambda^{n+1}_i&=\lambda_i^n-\Delta t(w_i^n-w_d)-\frac{\Delta t}{N}\sum_{j=1}^{N}R(w_i(t^n),w_j(t^n))\lambda^{n+1}_i,\;\;
\lambda_{i}^{n+1}=0.
\end{align}
\begin{align}\label{optimality}
&0= \Delta t\nu u^n+\frac{\Delta t}{N} \sum_{j=1}^N\lambda^{n+1}_j.
\end{align}
Upon substituting the terminal condition for $\lambda^{n+1}_j$ and expressing $u$ in terms of $\lambda^{n+1}_j$ we obtain 
the feedback control \eqref{eq:Dcontrol}.

\begin{remark}
In order to generalize the idea we may assume that the control acts differently on each agent. For example, one can consider the situation where action of the control $u$, acting on the single agents, is influenced  by the individual opinion. Therefore we replace $u$  in \eqref{eq:pbm} by $uQ(w_i)$, where $Q(w)$ is such that $q_m\leq Q(w)\leq q_M$.
Following the previous computation, the action of the control, at discrete time, is driven by
\begin{equation}
u^n Q_i^n = -\frac{\Delta t}{\nu N} \sum_{j=1}^N(w_j^{n+1}-w_d)Q^n_j Q_i^n
\label{eq:contag}
\end{equation}
where $Q^n_i=Q(w^n_i)$. Then the control dynamics on the opinion is  described by
\begin{equation}
\label{eq:DFwdQ}
\begin{aligned}
w^{n+1}_i&=w^n_i+\frac{\Delta t}{N}\sum_{j=1}^{N}P^n_{ij}(w^n_j-w^n_i)+\Delta t u^nQ_i^n,\qquad
w^0_i &= w_{0i}.
\end{aligned}
\end{equation}

\end{remark}

\section{Boltzmann description of constrained opinion consensus}
In this section, we consider a binary Boltzmann dynamic corresponding to the above model predictive control formulation. We emphasize that the assumption that opinions are formed mainly by binary interactions is rather common, see for example \cite{BS:09,GGL:12,PT:13,Tos06}.    Following \cite{AlbiPareschi2013ab, FoHaTo:11, PT:13} the first step is to reduce the dynamic to binary interactions.
Let consider the model predictive control system \eqref{eq:Dfwd}--\eqref{eq:Dcontrol} in the simplified case of only two interacting agents, numbered $i$ and $j$. Their opinions are modified in the following way
\begin{equation}\label{eq:DBin}
\begin{aligned}
&w^{n+1}_i=w_i^n+\frac{\Delta t}{2}P^n_{ij}(w_j^n-w^n_i)+\Delta t u^n,\\
&w^{n+1}_j=w_j^n+\frac{\Delta t}{2}P^n_{ji}(w_i^n-w^n_j)+\Delta t u^n,\\
\end{aligned}
\end{equation}
where the control
\begin{align}\label{eq:Icontrol}
&u^n=-\frac{\Delta t}{2\nu} \left((w^{n+1}_j-w_d)+(w^{n+1}_i-w_d))\right),
\end{align}
%\end{equation}
is implicitly defined in terms of the opinions pair at the time $n+1$. 
%can be interpreted as binary interaction between particles, in which the post-interaction opinions depends on the opinion of both interacting agents, moreover equation \eqref{eq:Icontrol} expresses the action of the control $u^n$, whit  implicit dependency  on the post-interaction opinions.  
The above linear system, however, can be easily inverted and its solutions can be written again in the form (\ref{eq:DBin}) where now the control is expressed explicitly in terms of the opinions pair at time $n$ as
\begin{align}\label{eq:Econtrol}
&u^n       = -\frac{1}{2}\frac{\Delta t}{\nu+\Delta t^2} \left((w^{n}_j-w_d)+(w^{n}_i-w_d))\right)-\frac{1}{2}\frac{\Delta t^2}{\nu+\Delta t^2} \left(P_{ij}-P_{ji}\right)(w_j^n-w_i^n).
\end{align}
%Where we gain an operator of the order of $\Delta t^2$.
Note that, as a result of the inversion of the $2\times 2$ matrix characterizing the linear system \eqref{eq:DBin}-\eqref{eq:Icontrol}, in the explicit formulation the control contains a term of order $\Delta t^2$.

\subsection{Binary interaction models}

In order to derive a kinetic equation we introduce a density distribution of particles  $f(w,t)$ depending on the opinion variable $w\in \I$ and time $t\geq0$.
The precise meaning of the density $f$ is the following. Given the
population of agents under study, if the opinions are defined on a subdomain $\Omega
\subset \I$ , the integral
 \[
 \int_\Omega f(w,t)\, dw
 \]
represents the number density of individuals with opinion included in $\Omega$ at time
$t> 0$. It is assumed that the density function is normalized to $1$, that
is
 \[
 \int_\I f(w,t) \, dw = 1.
 \]
The kinetic model can be derived by considering the change in time of $f(w,t)$ depending on the interactions with the other individuals. This change depends on the balance between the gain and loss due to the binary interactions. 

Accordingly to the explicit binary interaction \eqref{eq:DBin}, two agents with opinion $w$ and $v$ modify their opinion as 
\begin{equation}
\begin{aligned}\label{eq:BinD}
w^*&=\left(1-\alpha P(w,v)\right)w+\alpha P(w,v)v-\frac{\beta}{2}\left((v-w_d)+(w-w_d)\right)\\
&\qquad-\alpha\frac{\beta}{2}((P(w,v)-P(v,w))(w-v))+\Theta_1 D(w),\\
v^*&=\left(1-\alpha P(v,w)\right)v+\alpha P(v,w)w-\frac{\beta}{2} \left((v-w_d)+(w-w_d)\right)\\
&\qquad-\alpha\frac{\beta}{2}((P(v,w)-P(w,v))(v-w))+\Theta_2 D(w),\\
\end{aligned}
\end{equation}
where we included an additional noise term as in~\cite{Tos06}, to take into account effects falling outside the description of the model, like changes of opinion due to personal access to information. 
In \eqref{eq:BinD} we defined the following nonnegative quantities 
\begin{align}\label{eq:beta}
\alpha=\frac{\Delta t}{2},\qquad \beta = \frac{4\alpha^2}{\nu+4\alpha^2},
\end{align}
which represent the strength of the compromise and of the control respectively.
%
%The post-interaction opinions are determined by the following explicit binary relation
%\begin{equation}
%\begin{aligned}\label{eq:BinD}
%w^*&=w\left(1-\alpha\left( 1-\frac{\beta}{2}\right) P(w,v)-\frac{\beta}{2}\alpha P(v,w)\right)+v\left(\alpha\left(1-\frac{\beta}{2}\right)P(w,v)+\frac{\beta}{2}\alpha P(v,w)\right)\\&\qquad\qquad+\beta\left(w_d-\frac{w+v}{2}\right)+\Theta_1 D(w),\\
%v^*&=v\left(1-\alpha\left( 1-\frac{\beta}{2}\right) P(v,w)-\frac{\beta}{2}\alpha P(w,v)\right)+w\left(\alpha\left(1-\frac{\beta}{2}\right)P(v,w)+\frac{\beta}{2}\alpha P(w,v)\right)\\&\qquad\qquad+\beta\left(w_d-\frac{w+v}{2}\right)+\Theta_2 D(w).\\
%\end{aligned}
%\end{equation}
%%Where \begin{equation}\label{eq:beta}
%%\beta=\frac{4\alpha^2}{\nu+4\alpha^2},
%%\end{equation}
The noise term is characterized by the random variables $\Theta_1$ and $\Theta_2$ taking values on a set $\mathcal{B}\subset\mathbb{R}$, with identical distribution of mean  zero and variance $\sigma^2$ measuring the the degree of spreading of opinion due to diffusion. The function $D(\cdot)$ represents the local relevance of diffusion for a given opinion, and is such that $0\leq D(w)\leq 1$.

% Note that without any restriction on the parameters of \eqref{eq:BinD}, $w^*,v^*$ are not contained in the bounds, cause of forcing term to $w_d$ and the diffusion.
% In the next section we will skip this issue assuming directly that the bounds are always preserved by the binary interactions.
In the absence of diffusion, from \eqref{eq:BinD} it follows that 
\begin{subequations}
\begin{align}
w^*+v^* &= \left(1-\beta \right)(v+w)+2\beta w_d+\alpha(1-\beta)(P(w,v)-P(v,w))(v-w)\label{eq:meanop}\\
w^*-v^*   &=(w-v)(1-\alpha(P(w,v)+P(v,w))),\label{eq:conc}
\end{align}
\end{subequations}
thus in general the mean opinion is not conserved. Since $0\leq P(w,v)\leq 1$, if we assume $0\leq \alpha\leq 1/2$ from \eqref{eq:conc} we have
%\begin{equation}
\begin{align}
&|w^*-v^*|=\left(1-\alpha (P(w,v)+P(v,w)\right)|w-v|\leq(1-2\alpha)|w-v|,%\\\\
%&\epsilon(w,v)=\left(1-\alpha (P(w,v)+P(v,w))\right)\leq 1
\end{align}
which tells that the relative distance in opinion between two agents cannot increase after each interaction. 

When dealing with a kinetic problem in which the
variable belongs to a bounded domain we must deal with 
additional mathematical difficulties in the definition of
agents interactions. In fact, it is essential to consider only
interactions that do not produce values outside the finite
interval. The following proposition gives a sufficient condition to preserve the 
bounds. 

\begin{prop}
Let us assume that $0 < P(w,v)\leq 1$ and 
\begin{align}\label{eq:assumption}
 &\frac{\beta}{2}\leq\alpha p,\qquad|\Theta_i|< d\left(1-\frac{\beta}{2}\right),\quad i=1,2
\end{align}
where $p=\min_{w,v\in\I}\left\{P(w,v)\right\}>0$ and $d=\min_{w\in\I}\left\{(1-w)/D(w),D(w)\neq 0\right\}>0$, then the binary interaction \eqref{eq:BinD} preserves the bounds, i.e. the post-interaction opinions $w^*,v^*$ are contained in $\I=[-1,1]$.
\label{pr:1}
\end{prop}
%\noindent
\proof
%\textbf{Proof.} %We want to prove that relation \eqref{eq:BinD} keeps the bounds for $w^*$.
We will proceed in two subsequent steps, first by considering the case of interactions without noise and second by including the noise action.
Let us define the following quantity 
\begin{align}\label{eq:gamma}
\gamma = \alpha\left(1-\frac{\beta}{2}\right)P(w,v)+\alpha\frac{\beta}{2}P(v,w),
\end{align}
where $0\leq \beta \leq 1/2$ by definition.\\
Thus relation \eqref{eq:BinD} in absence of noise can be rewritten as 
\begin{align}
w^*=\left(1-\gamma-\frac{\beta}{2}\right)w+\left(\gamma-\frac{\beta}{2}\right)v+\beta w_d,
\end{align}
therefore it is sufficient that the following bounds are satisfied
 \begin{align}\label{eq:bounds}
 &\frac{\beta}{2}\leq\gamma\leq1-\frac{\beta}{2}
 \end{align}
to have a convex combination of $w$, $v$ and $w_d$.
From equation \eqref{eq:gamma}, by the assumption on $P(w,v)$, we have $\alpha p \leq \gamma \leq \alpha$. 
Therefore the left bound requires that $\alpha p\leq \beta/2$, which gives the first assumption in \eqref{eq:assumption}. 

If we now consider the presence of noise, we have
\begin{align}
\label{eq:bcon}
w^*=\left(1-\gamma-\frac{\beta}{2}\right)w+\left(\gamma-\frac{\beta}{2}\right)v+\beta w_d+D(w)\Theta_1.
\end{align}
%Without loss of generality because of symmetry, we consider $w\geq0$, then 
Equation \eqref{eq:bcon} implies the following inequalities
\begin{align*}
w^*\leq& \left(1-\gamma-\frac{\beta}{2}\right)w + \left(\gamma-\frac{\beta}{2}\right)+\beta w_d+D(w)\Theta_1\\
\leq& \left(1-\gamma-\frac{\beta}{2}\right)w + \left(\gamma+\frac{\beta}{2}\right)+D(w)\Theta_1.
\end{align*}
Finally, the last relation is bounded by one if 
$$\Theta_1\leq \left(1-\gamma-\frac{\beta}{2}\right)\frac{(1-w)}{D(w)}, \qquad D(w)\neq 0.$$
\noindent
which yields the second condition in \eqref{eq:assumption}. The same results are readily obtained for the post interacting opinion $v^*$.
\qed 
\begin{remark}
From the above proposition it is clear that agents should have a minimal amount of propensity to change their opinion in order for the control to act without risking to violate the opinion bounds. This reflects the fact that extreme opinions are very difficult to change and cannot be controlled in general without some additional assumption or model modification.
In the case of $\Theta_i= 0$, $\alpha \neq 0$ we obtain from \eqref{eq:assumption} the condition 
$$ \frac{2\alpha}{ \nu + 4\alpha^{2}}  \leq p. $$
This condition can be satisfied provided either $\alpha$ is sufficiently small or $\nu$ sufficiently large. 
%Note that for 
%large values of $\nu$ the control $u^{n}$ tends to zero according to \eqref{eq:pbic}. Hence, the condition on $\nu$ ensures
%that the bounds on $w$ are fulfilled.
\end{remark}

\subsection{Main properties of the Boltzmann description}
In general we can recover the time evolution of the density $f(w,t)$ through \eqref{eq:BinD} considering for a suitable test function $\varphi(w)$ an integro-differential equation of Boltzmann type in weak form \cite{PT:13}
\begin{align}\label{eq:Boltz}
\frac{d}{dt}\int_\I \varphi(w)f(w,t)dw=(Q(f,f),\varphi),
\end{align}
 where 
\begin{align}
\label{eq:Bcoll}
(Q(f,f),\varphi)=\left\langle\int_{\I^2}B_{int} \left(\varphi(w^*)-\varphi(w)\right)f(w,t)f(v,t)~dw~dv\right\rangle.
\end{align}
In \eqref{eq:Bcoll}, as usual, $\langle\, \cdot\, \rangle$ denotes the expectation with respect to the random variables $\Theta_i$, $i=1,2$ and the nonnegative interaction kernel $B_{int}$ is related to the probability of the microscopic interactions. The simplest choice which assures that the post interacting opinions preserves the bounds is given by  
\begin{align}
B_{int}=B_{int}(w,v,\Theta_1,\Theta_2)=\eta \chi(|w^*|\leq1)\chi(|v^*|\leq1)
\label{eq:ker}
\end{align}
where $\eta>0$ is a constant rate and $\chi(\,\cdot\,)$ is the indicator function.
A main simplification occurs if the bounds of $w^*,v^*$ are preserved by \eqref{eq:BinD} itself and the interaction kernel is independent on $w,v$, this will corresponds the classical Boltzmann equation for Maxwell molecules. In the rest of the paper, thanks to  Proposition \ref{pr:1}, we will pursue this direction. Following the derivation in \cite{CordierPareschiToscani2005aa,PT:13} the present results can be extended to kernels in the form \eqref{eq:ker}.

%\subsubsection{Main simplification and properties}

Let us assume that $|w^*|\leq1$ and $|v^*|\leq 1$, therefore the interaction dynamic of $f(w,t)$ can be described by the following Boltzmann operator
\begin{align}\label{eq:MBoltz}
(Q(f,f),\varphi)=\eta\left\langle\int_{\I^2} \left(\varphi(w^*)-\varphi(w)\right)f(w,t)f(v,t)~dw~dv\right\rangle.
\end{align}
%where $\eta$ represents a constant rate of interaction.
The above collisional operator guarantees the conservation of the total number of agents, corresponding to $\varphi(w)=1$, which is the only conserved quantity of the process. 
Let us remark that, since $f(w,t)$ is compactly supported in $\I$ then by conservation of the moment of order zero all the moments are bounded. By the same arguments in \cite{Tos06} the existence of a uniform bound on moments implies that the class of probability densities $\{f(w,t)\}_{t\geq 0}$
is tight, so that any sequence $\{f(w,t_n)\}_{t_n\geq 0}$ contains an infinite subsequence which
converges weakly as $t\to\infty$ to some probability measure $f_\infty$.

For $\varphi(w)=w$, we obtain the evolution of the average opinion. We have
\begin{align}\label{eq:Boltz1a}
\frac{d}{dt}\int_\I wf(w,t)dw=\eta\left\langle\int_{\I^2} \left(w^*-w\right)f(w,t)f(v,t)~dw~dv\right\rangle
\end{align}
or equivalently 
\begin{align}\label{eq:Boltz2a}
\frac{d}{dt}\int_\I wf(w,t)dw=\frac{\eta}{2}\left\langle\int_{\I^2} \left(w^*+v^*-w-v\right)f(w,t)f(v,t)~dw~dv\right\rangle.
\end{align}
Indicating the average opinion as  
\be
m(t)=\int_\I w f(w,t) \,dw,
\ee 
from relation \eqref{eq:Boltz2a} and \eqref{eq:meanop}, since $\Theta_i$, $i=1,2$ have zero mean, we obtain
\begin{align}
\nonumber
\frac{d}{dt}m(t)&=\frac{\eta}{2}\beta\int_{\I^2}\left(2w_d-w-v\right)f(v)f(w)~dw\,dv+\\
\nonumber
&\quad+\frac{\eta}{2}\alpha(1-\beta)\int_{\I^2}\left(P(w,v)-P(v,w)\right)(v-w)f(v)f(w)~dw\,dv\\
\label{eq:meanev}
    &=\eta\beta(w_d-m(t))+\eta\alpha(1-\beta)\int_{\I^2}(P(w,v)-P(v,w))vf(v)f(w)~dw\,dv.
\end{align}
Note that the above equation for a general $P$ is not closed.
Since $0\leq P(w,v)\leq 1$ we have $|P(w,v)-P(v,w)|\leq1$, then we can bound the derivative 
\begin{align*}
%\eta\beta(w_d-m)-\eta\alpham\leq& \frac{d}{dt}m\leq\eta\beta(w_d-m)+\eta\alpham\\
\eta\beta w_d-\eta(\beta+\alpha(1-\beta))m(t)\leq~ \frac{d}{dt}m(t)~\leq\eta\beta w_d-\eta(\beta-\alpha(1-\beta))m(t)
\end{align*}
solving on both sides we obtain the following estimate 
\begin{align*}
m(t)\geq\frac{\beta}{\beta+\alpha(1-\beta)}\left(1-e^{-\eta(\beta+\alpha(1-\beta))t}\right)w_d+m(0)e^{-\eta(\beta+\alpha(1-\beta))t}\\
m(t)\leq\frac{\beta}{\beta-\alpha(1-\beta)} \left(1-e^{-\eta(\beta-\alpha(1-\beta))t}\right)w_d+m(0)e^{-\eta(\beta-\alpha(1-\beta))t}.
\end{align*}
If we now assume that
\be
\nu<4\alpha,
\ee
then $\beta-\alpha(1-\beta)>0$ and if the average $m(t)\to m_{\infty}$ as $t\to\infty$ we have the bounds 
\be
\frac{4\alpha}{4\alpha+\nu}w_d \leq m_{\infty} \leq \frac{4\alpha}{4\alpha-\nu}w_d.
\ee
Therefore small values of $\nu$ force the mean opinion towards the desired state. 
%This gives a restriction on the choice of the regularization parameter $\nu$.
%\begin{remark}
In the symmetric case $P(v,w)=P(w,v)$, equation \eqref{eq:meanev} is in closed form and can be solved explicitly  
\begin{align}\label{eq:ave}
m(t)=\left(1-e^{-\eta\beta t}\right)w_d+m(0)e^{-\eta\beta t}
\end{align}
which in the limit $t\to\infty$ converges to $w_d$, for any choice of the control parameters. 

Let us now consider the case $\varphi(w)=w^2$ in the simplified situation of $P(w,v)=1$. We have
\begin{align}\label{eq:weaken}
\frac{d}{dt}\int_\I w^2 f(w,t)dw=\frac{\eta}{2}\left\langle\int_{\I^2} \left((w^*)^2+(v^*)^2-w^2-v^2\right)f(w,t)f(v,t)~dw~dv\right\rangle.
\end{align}
Denoting by
\be
E(t)=\int_\I w^2 f(w,t)dw,
\ee
easy computations show that
\begin{align}
\label{eq:ener}
\nonumber
\frac{d}{dt}E(t)=&-\eta\left(2\alpha(1-\alpha)+\beta\left(1-\frac{\beta}{2}\right) \right)(E(t)-m(t)^2)-2\eta\beta \left(\beta(m(t)^2-w_d^2)\right.\\[-.25cm]
\\[-.25cm]
\nonumber
&\left. +(1-\beta)m(t)(m(t)-w_d) \right)+\eta\sigma^2\int_\I D(w) f(w,t)\,dw, 
\end{align}
where we used the fact that $\Theta_i$, $i=1,2$ have zero mean and variance $\sigma^2$. In absence of diffusion, since $m(t)\to w_d$ as $t\to\infty$, we obtain that $E(t)$ converges exponentially to $w_d^2$ for large times. Therefore the quantity
\be
\int_\I f(w,t) (w-w_d)^2\,dv=E(t)^2+w_d^2-2m(t)w_d,
\ee
goes to zero as $t \to\infty$. This shows that, under the above assumptions, the steady state solution has the form of a Dirac delta $f_\infty(w)=\delta(w-w_d)$ centered in the desired opinion state.  

%\end{remark}

\section{Fokker-Planck modeling}
In the general case,  it is quite difficult to obtain analytic results on the large time behavior of the kinetic equation \eqref{eq:MBoltz}. 
As it is usual in kinetic theory, particular asymptotic limit of the Boltzmann model result in simplified models, generally of Fokker-Planck type, for which the study of the theoretical properties is often easier \cite{PT:13}. 
%Here we consider first the so-called quasi invariant opinion limit of the Boltzmann model and then a mean field limit of the microscopic dynamic.
\subsection{The quasi-invariant opinion limit}
\label{qiol}
%In order to recover a simpler operator then \eqref{eq:MBoltz}, we look to asymptotic properties of the kinetic equation. 
The main idea is to rescale the interaction frequency $\eta$, the propensity strength $\alpha$, the diffusion variance $\sigma^2$ and the action of the control $\nu$ at the same time, in order to maintain at level of the asymptotic procedure the memory of the microscopic interactions \eqref{eq:BinD}. This approach is usually referred to as {quasi--invariant opinion} limit~\cite{PT:13,Tos06} and is closely related to the {grazing collision limit} of the Boltzmann equation for Coulombian interactions (see \cite{FPTT12, Vill98}). 

We make the following scaling assumptions 
\begin{align}\label{eq:scale}% scaling
\alpha=\epsi,\qquad
\eta = \frac{1}{\epsi},\qquad
\sigma^2=\epsi\varsigma, \qquad \nu =\epsi\kappa,
\end{align}
where $\epsi>0$ and as a consequence the coefficient $\beta$ in \eqref{eq:BinD} takes the form
\begin{align*}
\beta = \frac{4\epsi}{\kappa +4\epsi}.
\end{align*}
This corresponds to the situation where the interaction operator concentrates on binary interactions which produce a very small change in the opinion of the agents. From a modeling viewpoint, we require that scaling \eqref{eq:scale} in the limit $\varepsilon\to 0$ preserves the main macroscopic properties of the kinetic system. To this aim, let us observe that the evolution of the scaled first two moments for $P(w,v)=1$ reads
\begin{align*}
\frac{d}{dt}m(t)=&\frac{4}{\kappa+4\varepsilon}(w_d-m(t)),\\
\nonumber
\frac{d}{dt}E(t)=&-2 \left((1-\varepsilon)+\frac{2}{\kappa+4\varepsilon}\left(1-\frac{2\varepsilon}{\kappa+4\varepsilon}\right) \right)(E(t)-m(t)^2)\\
&-\frac{8}{\kappa+4\varepsilon} \left(\frac{4\varepsilon}{\kappa+4\varepsilon}(m(t)^2-w_d^2)
+\left(1-\frac{4\varepsilon}{\kappa+4\varepsilon}\right)m(t)(m(t)-w_d) \right)\\
\nonumber
&+\varsigma\int_\I D(w) f(w,t)\,dw, 
\end{align*}
which in the limit $\varepsilon\to 0$ gives
\begin{align}
\label{eq:ms1}
\frac{d}{dt}m(t)=&\frac{4}{\kappa}(w_d-m(t)),\\
\nonumber
\frac{d}{dt}E(t)=&-2\left(1+\frac{2}{\kappa}\right)(E(t)-m(t)^2)\\[-.25cm]
\label{eq:es1}
\\[-.25cm]
\nonumber
&-\frac{8}{\kappa}m(t)(m(t)-w_d)+\varsigma\int_\I D(w) f(w,t)\,dw. 
\end{align}
This shows that in order to keep the effects of the control and the diffusion in the limit it is essential that both $\nu$ and $\sigma^2$ scale as $\varepsilon$.

%In view of the above definitions this corresponds also to scale ${\Delta t}$ as $\varepsilon$. 
%and similarly the other model parameters in order to preserve the model behavior as $\epsilon\to 0$. This can be better understood by observing that for symmetric $P$ the scaled mean opinion satisfies the equation 
%Note that, as result of the scaling \eqref{eq:scale}, conditions \eqref{eq:assumption} are automatically satisfied in the limit $\varepsilon\to 0$, since both $\alpha\to 0$ and $\beta\to 0$, and therefore no lower bound on $P$ is necessary.
%%Through this approach we want to maintain at level of the asymptotic the memory of the microscopic interactions.
In the sequel we show how this approach leads to a constrained Fokker--Planck equation for the description of the opinion distribution. Even if our computations are formal, following the same arguments in \cite{PT:13,Tos06} it is possible to give a rigorous mathematical basis to the derivation. Here we omit the details for brevity.

\noindent
The scaled equation  \eqref{eq:MBoltz} reads
\begin{align}\label{eq:MBoltzs}
\frac{d}{dt}\int_\I \varphi(w)f(w,t)dw=\frac{1}{\epsi}\left\langle\int_{\I^2} \left(\varphi(w^*)-\varphi(w)\right)f(w,t)f(v,t)~dw~dv\right\rangle
\end{align}
where the scaled binary interaction dynamic \eqref{eq:BinD} can be written as
\begin{equation}\label{eq:BinScale}
w^* -w=\epsi P(w,v)(v-w)  + \frac{2\epsi}{\kappa+4\epsi}\left( 2w_d -(w+v)\right)+\Theta^{\varepsilon}_1 D(w)+O(\epsi^2),\\
\end{equation}
\noindent
where $\Theta^{\varepsilon}_1$ is a random variable with zero mean and variance $\varepsilon\varsigma$.

In order to recover the limit as $\epsi\to 0$ we consider the second-order Taylor expansion of $\varphi$ around $w$
\begin{align}\label{eq:Tay2}
\varphi(w^*)-\varphi(w)=(w^*-w)\varphi'(w)+\frac{1}{2}(w^*-w)^2\varphi''(\tilde{w})
\end{align}
where for some $0\leq\vartheta\leq1$ ,
$$\tilde{w}=\vartheta w^*+(1-\vartheta)w.$$ 
Therefore, inserting this expansion in the interaction integral \eqref{eq:MBoltzs} we get
\begin{equation}\label{eq:limBoltz}
\begin{aligned}
&\frac{1}{\epsi}
\left\langle\int_{\I^2} \left(\left(w^*-w\right)\varphi'(w)+\frac{1}{2}\left(w^*-w\right)^2\varphi''(w) \right)f(w)f(v)~dwdv\right\rangle+R(\epsi).
\end{aligned}
\end{equation}
The term $R(\epsi)$  denotes the  remainder and is given by
\begin{equation}
R(\epsi)=\frac{1}{2\epsi}\left\langle\int_{\I^2}\left(w^*-w\right)^2(\varphi''(\tilde{w})-\varphi''(w))f(w)f(v)~dwdv\right\rangle.
\label{eq:rest}
\end{equation}
Using now \eqref{eq:BinScale} we can write
\begin{equation}\label{eq:limBoltz2}
\begin{aligned}
&\frac{1}{\epsi}
\int_{\I^2} \left[\left(P(w,v)(v-w)+\frac{2\epsi}{\kappa+4\epsi}\left( 2w_d -(w+v)\right) \right)\varphi'(w)\right.\\
&\qquad\left.+\frac{\varsigma}{2}D(w)^2 \varphi''(w)\right]f(w)f(v)~dwdv+R(\epsi)+O(\varepsilon),
\end{aligned}
\end{equation}
where we used the fact that $\Theta_1^{\varepsilon}$ has zero mean and variance $\varepsilon\varsigma$.

By the same arguments in \cite{Tos06} it is possible to show rigorously that \eqref{eq:rest} converges to zero as soon as $\epsi\to 0$. 
Therefore we have as limiting operator of \eqref{eq:MBoltz} the following
\begin{align*}
\frac{d}{dt}\int_\I \varphi(w)f(w,t)dw=&\int_{\I^2}\left(P(w,v)(v-w)  + \frac{4}{\kappa}\left( w_d -\frac{w+v}{2}\right)\right)\varphi'(w)f(w)f(v)~dwdv\\
%\\&\qquad
&+\frac{\varsigma}{2}\int_{\I} D(w)^2\varphi''(w)f(w)~dw.
\end{align*} 
Integrating back by parts the last expression we obtain the following Fokker--Planck equation 
\begin{equation}
\begin{aligned}\label{eq:FP1}
\frac{\partial}{\partial t} f &+\frac{\partial}{\partial w}\mathcal{H}[f](w)f(w)+ \frac{\partial}{\partial w}\mathcal{K}[f](w)f(w)~dv=\frac{\varsigma}{2}\frac{\partial^2}{\partial w^2}(D(w)^2f(w)),
\end{aligned}
\end{equation}
where
\begin{align}
&\mathcal{K}[f](w)=\int_\I P(w,v)(v-w)f(v)~dv,\\
&\mathcal{H}[f](w)=\frac{4}{\kappa}\int_\I \left(w_d-\frac{w+v}{2}\right)f(v)~dv=\frac{4}{\kappa}\left(w_d-\frac{w+m}{2}\right).
\end{align}
%where $m=\int_\I wf(w)dw$.

\begin{remark}
%Note that 
%the scale assumption  \eqref{eq:scale}, is equivalent to  consider a scale of the time variable
%\begin{equation}\label{eq:timescale}
%\tau=\epsi t, \qquad g(w,\tau)=f(w,t)
%\end{equation}
%which implies $f_0(w)=g_0(w)$.  
The ratio between $\sigma^2/\alpha=\varsigma$ is of paramount importance in order to obtain in the limit the contribution of both controlled compromise propensity and  diffusion \cite{Tos06}.
Other limiting behaviors can be considered like diffusion dominated $(\varsigma\to \infty)$ or controlled compromise dominated $(\varsigma\to 0)$. 
%Moreover, since  $t=\tau/\epsi$ for $\epsi \to 0$, the scaled density $g(w,\tau)$ describes the large time behavior of $f(w,t)$.
\end{remark}
\begin{figure}[t]
\centering
\includegraphics[scale=0.4]{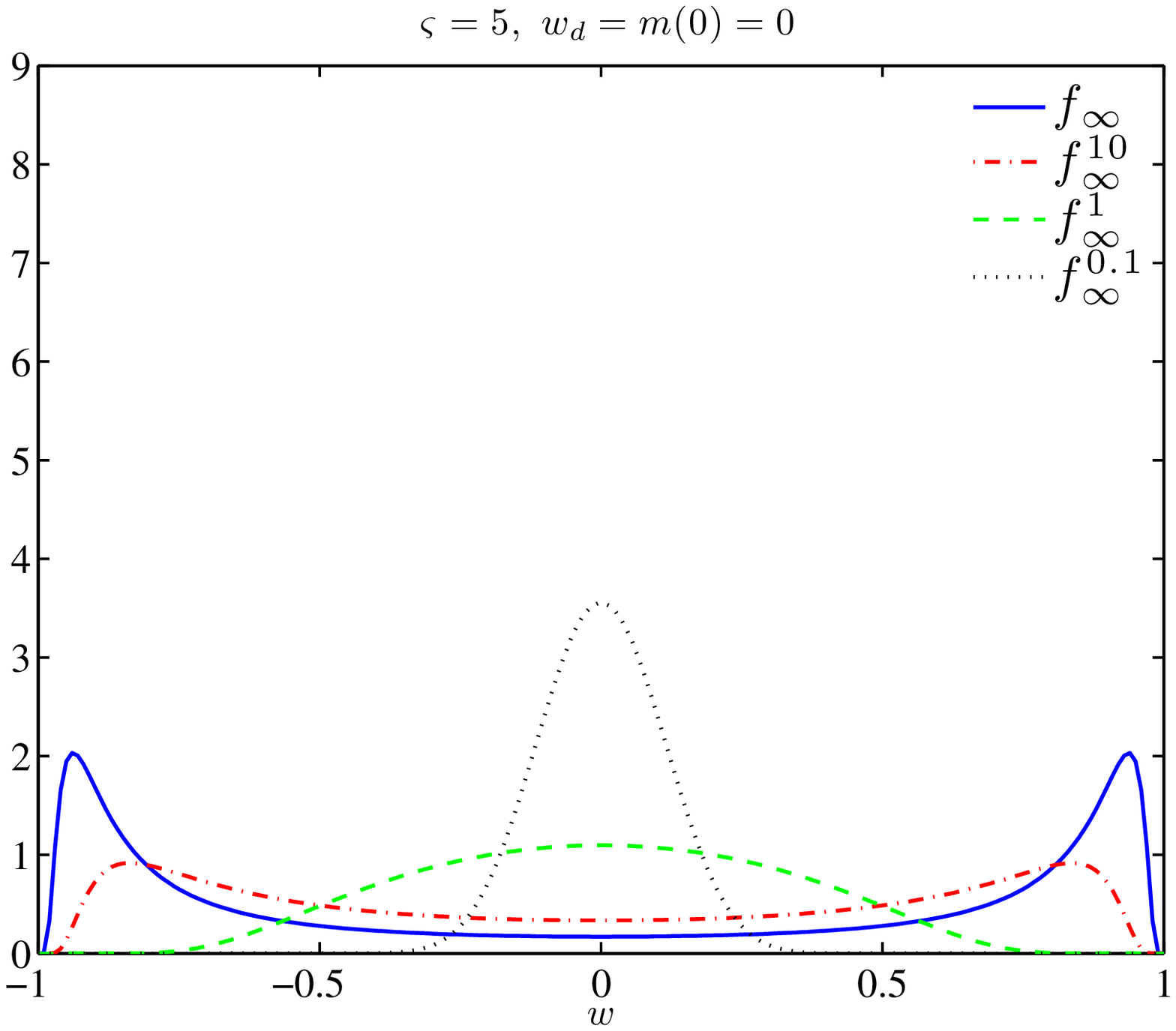}
%&
\includegraphics[scale=0.4]{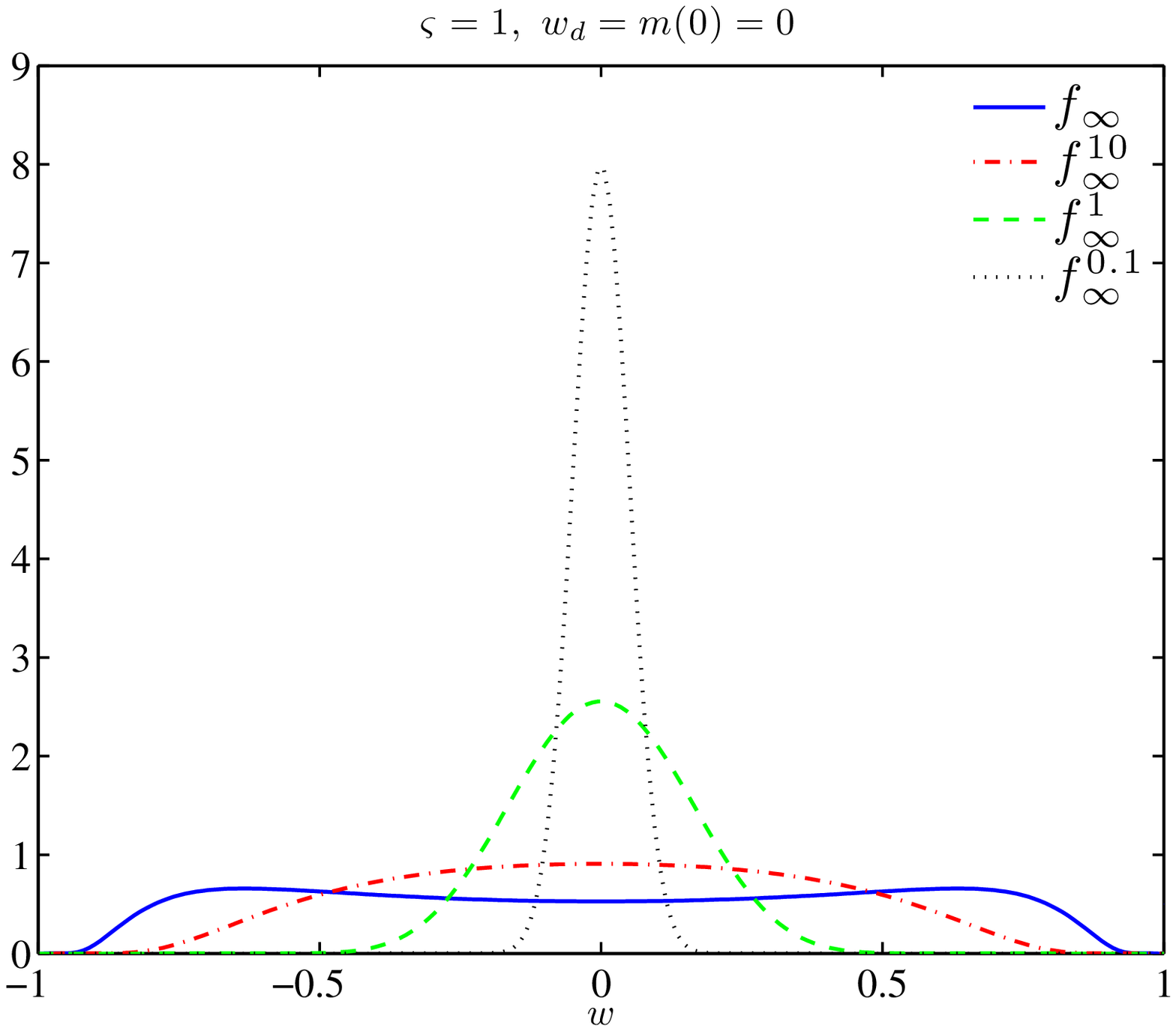}
\caption{Continuous line and dashed lines represent  the steady solutions $f_\infty$ and $f_\infty^\kappa$, respectively.
 On the left $w_d=m(0)=0$ with diffusion parameter $\varsigma=5$, on the right  $w_d=m(0)=0$ with diffusion parameter $\varsigma=2$.
In both cases the steady solution changes from a bimodal distribution to an unimodal distribution around $w_d$.}\label{fig:F1}
\end{figure}
\subsection{Stationary solutions}\label{sec:Steady}
In this section we analyze the steady solutions of the Fokker--Planck model \eqref{eq:FP1}, for particular choices of the microscopic interaction of the Boltzmann dynamic.
% for which is possible to derive the Fokker--Planck model and compute explicitly the steady state. 

Let consider the case in which $P(w,v)=1$. 
%In absence of control, i.e. $\kappa\to\infty$, the mean opinion is conserved $m=m(0)$ and steady solutions of \eqref{eq:FP1} satisfy the following differential equation, see \cite{Tos06} 
%\begin{align}\label{eq:edo}
%\frac{\varsigma}{2}\partial_w(D(w)^2f)=\left(m-w\right)f.
%\end{align}
In presence of the control the average opinion in general is not conserved in time, but since $m(t)$ converges exponentially in time to $w_d$, the steady state opinion solves
\begin{align}\label{eq:edo2}
&\frac{\varsigma}{2}\partial_w(D(w)^2f)=\left(1+\frac{2}{\kappa}\right)(w_d-w)f.
\end{align}
%  then the average opinion $m(t)$ evolves according to 
%\begin{align}\label{eq:FPave}
%m(t)=\left(1-e^{-4/\kappa t}\right)w_d+e^{-4/\kappa t}m(0),
%\end{align}
%which is obtained from the scaled equation \eqref{eq:Boltz2} through the quasi-invariant opinion limit.
If we now consider as diffusion function $D(w)=(1-w^2)$, then it is possible explicitly compute the solution of \eqref{eq:edo2} as follows \cite{Tos06}
\begin{align}\label{eq:sol}
f^{\kappa}_\infty(w)=&~\frac{C_{w_d,\varsigma,\kappa}}{(1-w^2)^2}\left(\frac{1+w}{1-w}\right)^{m/(2\varsigma)}\exp\left\{-\frac{1-m w}{\varsigma\left(1-w^2\right)}\left(1+\frac{2}{\kappa}\right)\right\}
\end{align}
where $C_{w_d,\varsigma,\kappa}$ is a normalization constant such that $\int f_{\infty}\,dw=1$. We remark that the solution is such that $f(\pm1)=0$, moreover due to the general non symmetry of $f$, the desired state reflects on the steady state through the mean opinion. Note that in the case $\kappa\to\infty$ we obtain the steady state of the uncontrolled equation \cite{Tos06}. We denote by $f_{\infty}(w)$ this latter uncontrolled stationary behavior. 
\begin{figure}[t]
\centering
%\includegraphics[scale=0.4]{Stationary_wd_k1.eps}
%&
\includegraphics[scale=0.4]{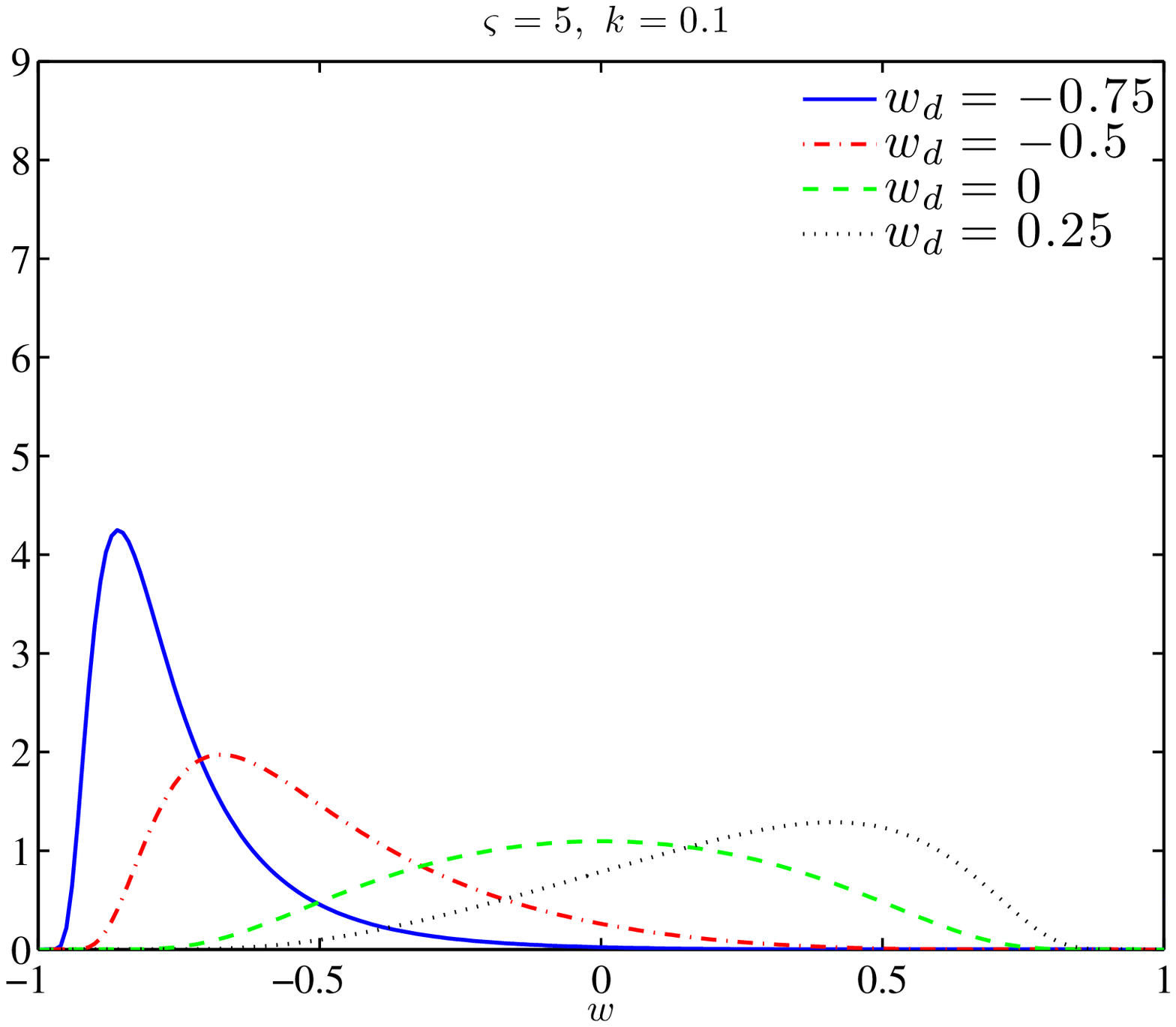}
%&
%\\
%\includegraphics[scale=0.4]{Stationary_wd_k3.eps}
%&
\includegraphics[scale=0.4]{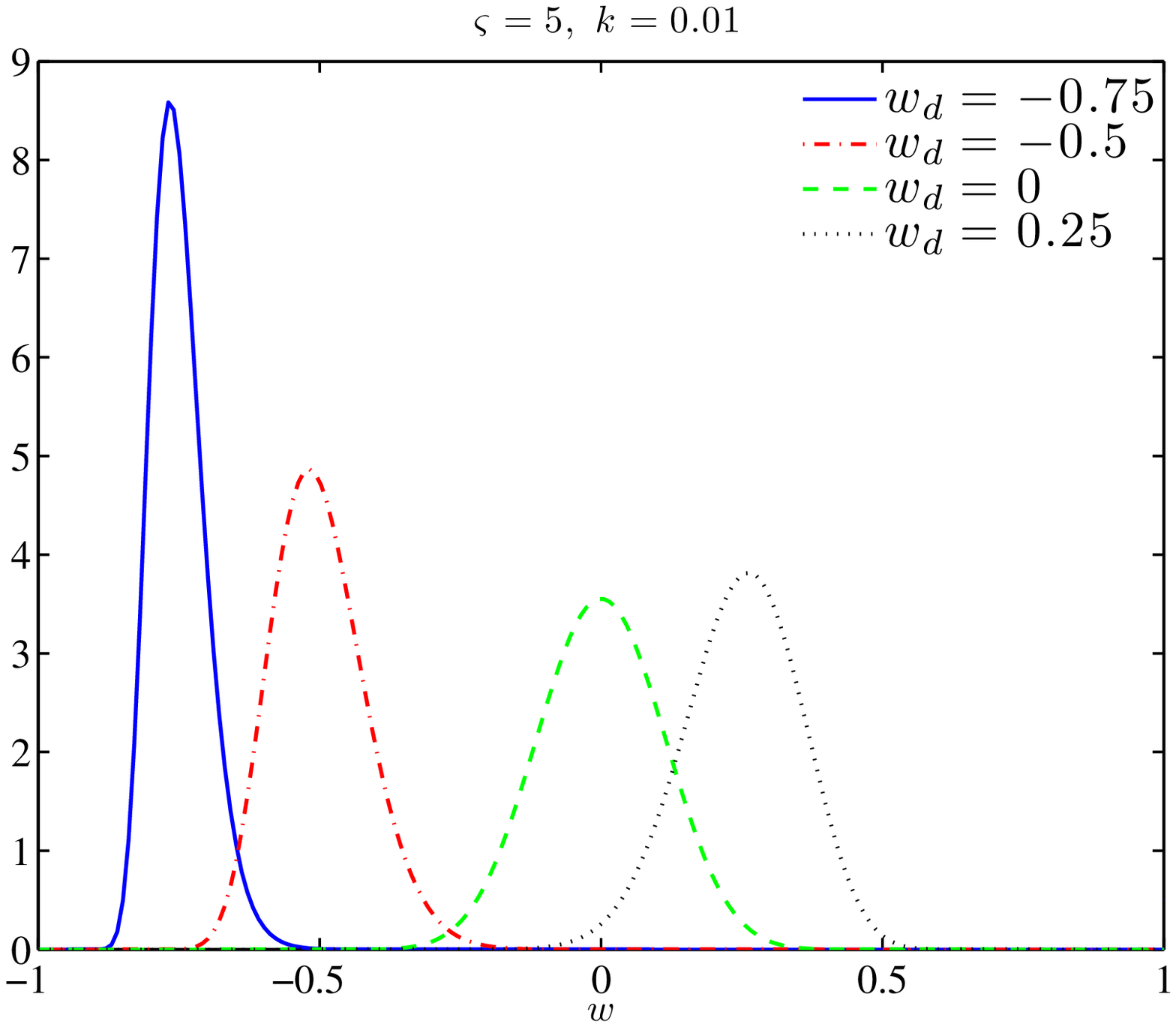}
\caption{Steady state solutions in the controlled case for different values of $\kappa$ and $w_d$. From left to right we change values of $\kappa=0.1$ and $\kappa=0.01$ for a fixed value of $\varsigma=5$ and different desired states $w_d=\left\{-0.75, -0.5, 0, 0.25\right\}$.}\label{fig:F2}
\end{figure}
%The explicit solution of \eqref{eq:edo2} can be generalized out of the solution of \eqref{eq:edo} as
%\begin{align}\label{eq:sol1}
%f_\infty^\kappa(w)& = \frac{C_{w_d,\varsigma,\kappa}}{(1-%w^2)^2}\left(S_{w_d,\varsigma}(w)\right)^{1+2/\kappa}
%\end{align}
%where $C_{w_d,\varsigma,\kappa}$ is  the normalization constant.
We plot in Figure \ref{fig:F1} the steady profile $f_\infty$ and $f_\infty^\kappa$ for different choices of the parameters $\kappa$ and $\varsigma$. The initial average opinion $m(0)$ is taken equal to the desired opinion $w_d$, in this way we can see that for $\kappa\to \infty$ the constrained steady profile approaches the unconstrained one, $f^\kappa_\infty\to f_\infty$. On the other hand small values of $\kappa$ give the desired distribution concentrated around $w_d$.
%In general  for steady solution we can not switch from $f_\infty$ to $f_\infty^\kappa$ playing with parameter $\kappa$, since the memory on the initial average opinion is lost for any $\kappa>0$.

In Figure \ref{fig:F2} we show the steady profile $f_\infty^\kappa$ for different choice of the parameters $\kappa$ and the desired state $w_d$. We can see that decreasing the value of $\kappa$ lead the profiles to concentrate around the requested value of $w_d$. 
\begin{figure}[t]
\centering
\includegraphics[scale=0.4]{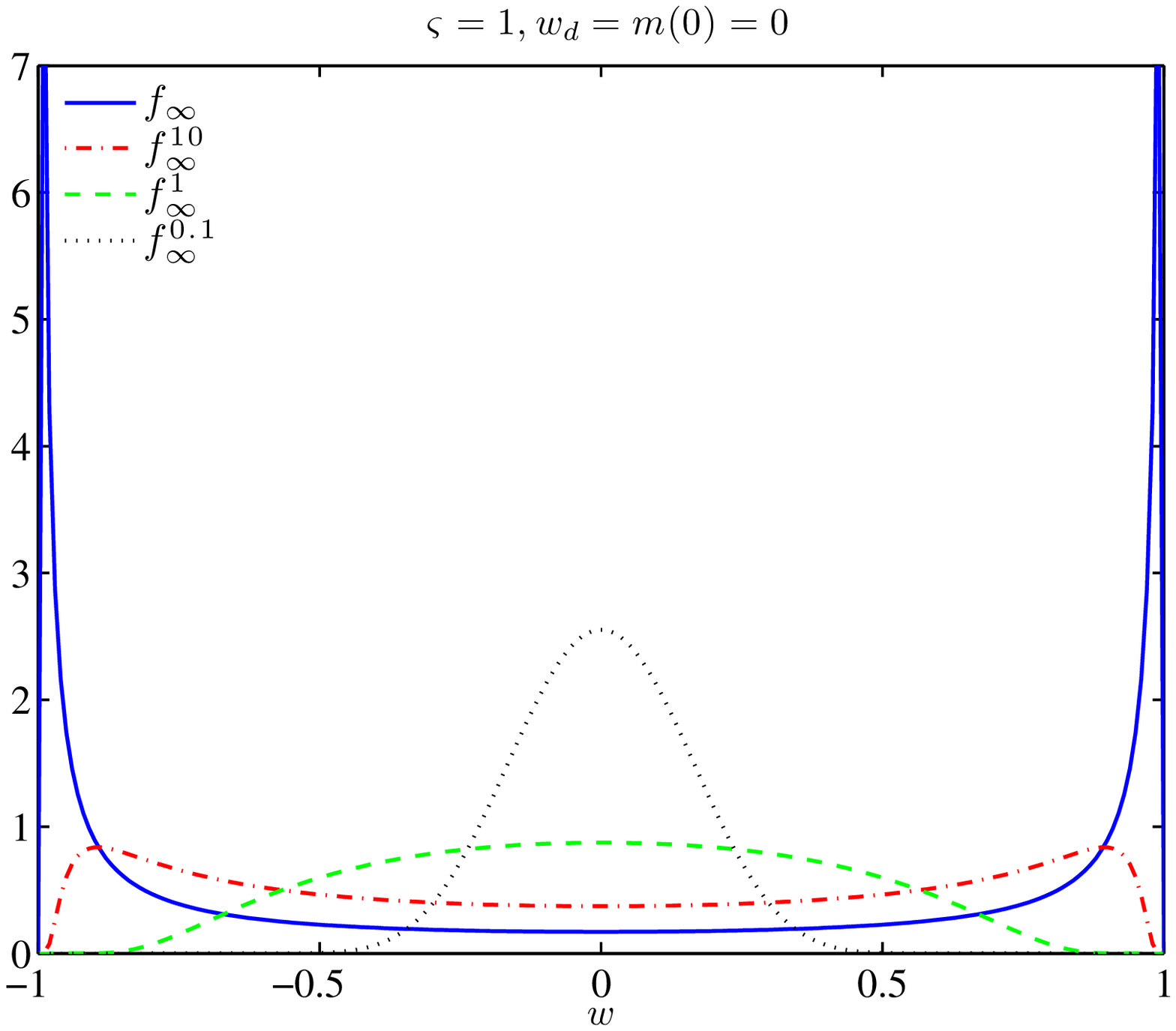}
%&
\includegraphics[scale=0.4]{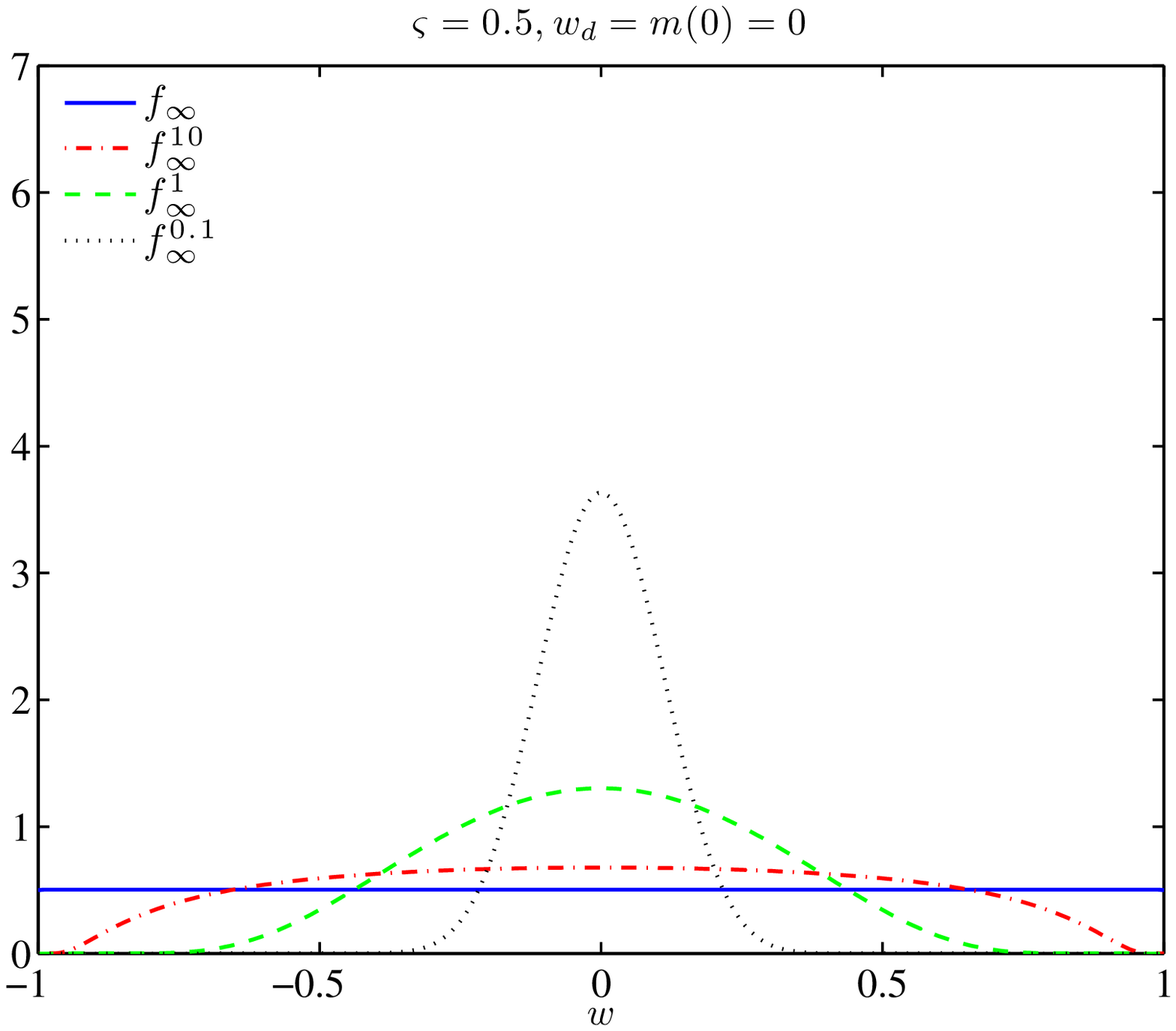}
\caption{Continuous line and dashed lines represent  the steady solutions $f_\infty$ and $f_\infty^\kappa$, respectively. On the left $w_d=m(0)=0$ with diffusion parameter $\varsigma=0.9$, on the right  $w_d=m(0)=0$ with diffusion parameter $\varsigma=0.5$, in this last case note that $f_\infty$ is a uniform distribution on $[-1,1]$.}\label{fig:F3}
\end{figure}

Let consider $P(w,v)=P(w)$ then stationary solutions of \eqref{eq:FP1} satisfy the following 
%\begin{align}\label{eq:edo3}
%&\frac{\varsigma}{2}\partial_w(D(w)^2f)=P(w)(m(0)-w)f.
%\end{align}
%In presence of control we have 
\begin{align}\label{eq:edo4}
&\frac{\varsigma}{2}\partial_w(D(w)^2f)=\left(P(w)+\frac{2}{\kappa}\right)(w_d-w)f.
\end{align}
%Second example we consider $D(w)=\sqrt{1-w^2}$ [...]
Taking $P(w)=1-w^2$ and $D(w)=1-w^2$ we can compute \cite{Tos06}
%\begin{align}
%f_\infty(w)=C_{\varsigma,m}(1-w)^{-2-\frac{m(0)-1}{\varsigma}}%(1+w)^{-2+\frac{m(0)+1}{\varsigma}}
%\end{align}
%with control
\begin{align}
f_\infty^\kappa(w) = C_{\varsigma,m}(1-w)^{-2-\frac{w_d-1}{\varsigma}-\frac{w_d}{\kappa\varsigma}}(1+w)^{-2+\frac{w_d+1}{\varsigma}+\frac{w_d}{\kappa\varsigma}}\exp\left\{-\frac{2}{\kappa}\frac{1-w_d w}{\varsigma\left(1-w^2\right)}\right\}
\label{eq:steady2}
\end{align}
We present in Figure \ref{fig:F3} different profiles of  $f^\kappa_\infty$ for $m(0)=w_d$, where we switch from the steady profile of the uncontrolled case to the steady profile \eqref{eq:steady2}.

\section{Other constrained kinetic models}
The constrained binary collision rule \eqref{eq:BinD} admits several variants accordingly to the different ways we realize the diffusion and control dynamics. 

From the modeling point of view we decided to introduce noise at the level of the explicit binary formulation \eqref{eq:DBin},\eqref{eq:Econtrol} as an external factor which can not be affected by the opinion maker. 
%Of course, the same could have been done directly at the level of the model predictive microscopic model \eqref{eq:Dfwd}-\eqref{eq:Dcontrol} after solving it explicitly as a linear system to obtain $w_i^{n+1}$.
In contrast, adding noise from the very beginning in \eqref{eq:pbm}-\eqref{eq:pbc}, or equivalently in the implicit formulation \eqref{eq:Dfwd}-\eqref{eq:Dcontrol}, would imply a different action of the control over the spreading of the noise. More precisely, for the binary interaction model this will originate the dynamic 
\begin{equation}
\begin{aligned}\label{eq:BinD2}
w^*&=\left(1-\alpha P(w,v)\right)w+\alpha P(w,v)v-\frac{\beta}{2}\left((v-w_d)+(w-w_d)\right)\\
&\qquad-\alpha\frac{\beta}{2}((P(w,v)-P(v,w))(w-v))+\left(1-\frac{\beta}{2}\right)\Theta_1 D(w)-\frac{\beta}{2}\Theta_2 D(v),\\
v^*&=\left(1-\alpha P(v,w)\right)v+\alpha P(v,w)w-\frac{\beta}{2} \left((v-w_d)+(w-w_d)\right)\\
&\qquad-\alpha\frac{\beta}{2}((P(v,w)-P(w,v))(v-w))+\left(1-\frac{\beta}{2}\right)\Theta_2 D(v)-\frac{\beta}{2}\Theta_1 D(w).\\
\end{aligned}
\end{equation}
For this binary dynamic preservation of the bounds is more delicate and the corresponding Boltzmann model is typically written using the kernel \eqref{eq:ker}. Note, however, that in the quasi-invariant opinion limit due to the rescaling \eqref{eq:scale} we have $\beta\to 0$ and therefore the limiting Fokker-Planck equation is again \eqref{eq:FP1}. 

Next we remark that the microscopic constrained system \eqref{eq:Dfwd}-\eqref{eq:Dcontrol} can be written in explicit form by solving the corresponding linear system for $w_1^{n+1},\ldots,w_N^{n+1}$. Straightforward computations yields the explicit formulation
\begin{equation}
w^{n+1}_i=w^n_i+\frac{\Delta t}{N}\sum_{j=1}^{N}P^n_{ij}(w^n_j-w^n_i)+\Delta t u^n,\qquad w^0_i= w_{0i},\label{eq:Dfwd2}
\end{equation}
where now
\begin{eqnarray}
u^n = \frac{(\Delta t)^2}{\nu + (\Delta t)^2}  \left(\frac1{N^2}\sum_{h,j=1}^N P(w_h,w_j)(w_j^n-w_h^n) \right)+\frac{\Delta t}{\nu + (\Delta t)^2} (w_d - m^n),
\label{eq:Econtrol2}
\end{eqnarray}
and we denoted by \[ m^n=\frac1{N}\sum_{j=1}^N w_j^n\] the mean opinion value. This show that a different way to realize the constrained binary dynamic \eqref{eq:BinD} is given by 
\begin{equation}
\begin{aligned}\label{eq:BinD3}
w^*&=\left(1-\alpha P(w,v)\right)w+\alpha P(w,v)v-\beta\left(m(t)-w_d\right)\\
&\qquad-\alpha\frac{\beta}{2}((P(w,v)-P(v,w))(w-v))+\Theta_1 D(w),\\
v^*&=\left(1-\alpha P(v,w)\right)v+\alpha P(v,w)w-\beta\left(m(t)-w_d\right)\\
&\qquad-\alpha\frac{\beta}{2}((P(v,w)-P(w,v))(v-w))+\Theta_2 D(w).\\
\end{aligned}
\end{equation}
Again preservation of the bounds is a difficult task and the Boltzmann equation is written in the general form \eqref{eq:Bcoll}. Performing the same computations as in Section \ref{qiol} we obtain the limiting Fokker-Planck equation \eqref{eq:FP1} with the simplified control term
\begin{align}
\mathcal{H}[f](w)=\frac{4}{\kappa}\left(w_d-m\right).
\end{align} 
The main difference now, is that when $m(t)\to w_d$ the contribution of the control vanish, $\mathcal{H}[f](w)\to 0$, and the steady states corresponds to those of the unconstrained equation by Toscani \cite{Tos06} in the case where the mean opinion is given by the desired state. In other words, in the examples of Section \ref{sec:Steady}, they are given by \eqref{eq:sol} and \eqref{eq:steady2} in the limit case $\kappa\to\infty$. Therefore, in this case, the action of the control is weaker, since it is not able to act on any opinion distribution with mean opinion given by the desired state.

Finally, from system \eqref{eq:contag}-\eqref{eq:DFwdQ}, we can also generalize \eqref{eq:BinD} with an agent dependent action of the control. Following the same derivation as in Section 3 we have the binary interaction rule
\begin{equation}
\begin{aligned}\label{eq:BinD4}
w^*&=\left(1-\alpha P(w,v)\right)w+\alpha P(w,v)v-\frac{\beta(w,v)}{2}\left(Q(v)(v-w_d)+Q(w)(w-w_d)\right)\\
&\qquad-\alpha\frac{\beta(w,v)}{2}(Q(w)P(w,v)-Q(v)P(v,w))(v-w)+\Theta_1 D(w),\\
v^*&=\left(1-\alpha P(v,w)\right)v+\alpha P(v,w)w-\frac{\beta(v,w)}{2}\left(Q(v)(v-w_d)+Q(w)(w-w_d)\right)\\
&\qquad-\alpha\frac{\beta(v,w)}{2}(Q(v)P(v,w)-Q(w)P(w,v))(w-v)+\Theta_2 D(v),\\
\end{aligned}
\end{equation}
where
\[
\beta(w,v)=\frac{4\alpha^2Q(w)}{\nu+2\alpha^2(Q(v)^2+Q(w)^2)},
\]
with property $\beta(w,v)Q(v)=\beta(v,w)Q(w)$. In this case, sufficient condition for the preservation of the bounds can be found provided that a minimal action of the control is admitted by the agents, namely assuming that $0<Q(\,\cdot\,)\leq 1$. Under the scaling \eqref{eq:scale} we obtain the general Fokker-Plank equation \eqref{eq:FP1} where now the control term reads
\be
\mathcal{H}[f](w)=\left(\frac{2}{\kappa}\int_\I\left(Q(w)(w_d-w)+Q(v)(w_d-v)\right)f(v)\,dv\right)Q(w).
\ee

%, it is in good agreement with the original model predictive control approximation (\ref{eq:Dfwd2})-(\ref{eq:Econtrol2}). 

%In fact, the steady states 
%satisfy
%\be
%\frac{\varsigma}{2}\partial_w(D(w)^2f)=(w_d-w)f,
%\label{eq:ss}
%\ee
%since again we have that $m(t)$ satisfies \eqref{eq:ms1}.

\section{Numerical examples}
In this section we report some numerical test obtained by solving the constrained Boltzmann equation with the binary interaction rule \eqref{eq:BinD} for different kind of opinion models. In the numerical simulations we use a Monte Carlo methods as described in Chapter 4 of \cite{PT:13}. 
We simulate equation \eqref{eq:FP1} for particular choices of the parameters of the model comparing the stationary solutions obtained in absence of control \cite{Tos06, ANT07} with different increasing actions of the control term. 

\subsubsection*{Quasi-invariant opinion limit}
In the first numerical example we compare the solutions obtained with the Monte Carlo method in the quasi-invariant opinion limit with the exact profile of the steady solution of the Fokker--Planck model \eqref{eq:FP1}. We consider the particular case 
\be
P(w,v)=1,\qquad D(w)=1-w^2, 
\ee
then exact solutions are described by \eqref{eq:sol}.
\begin{figure}[t]
\centering
\includegraphics[scale=0.35]{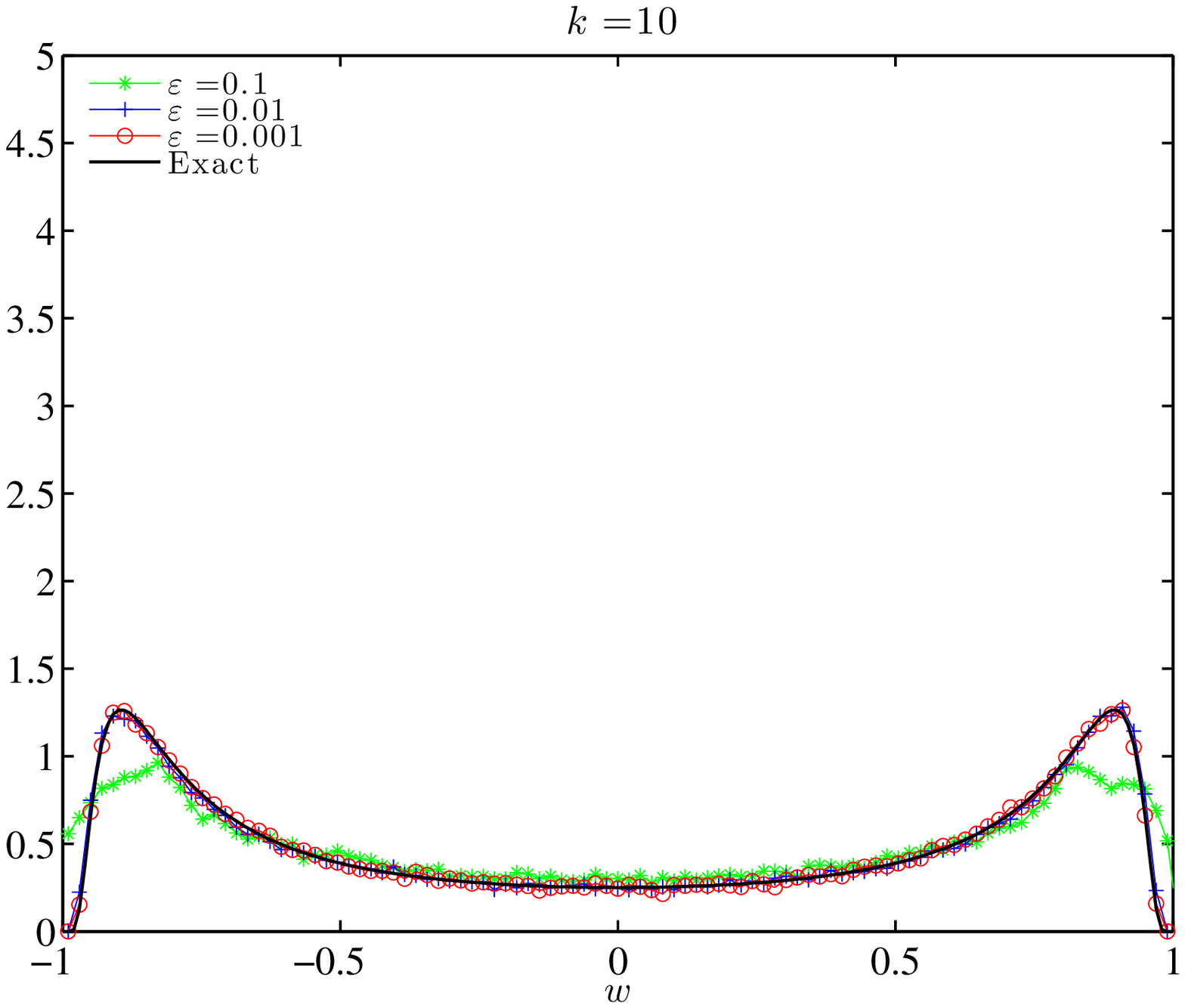}
%&
\includegraphics[scale=0.35]{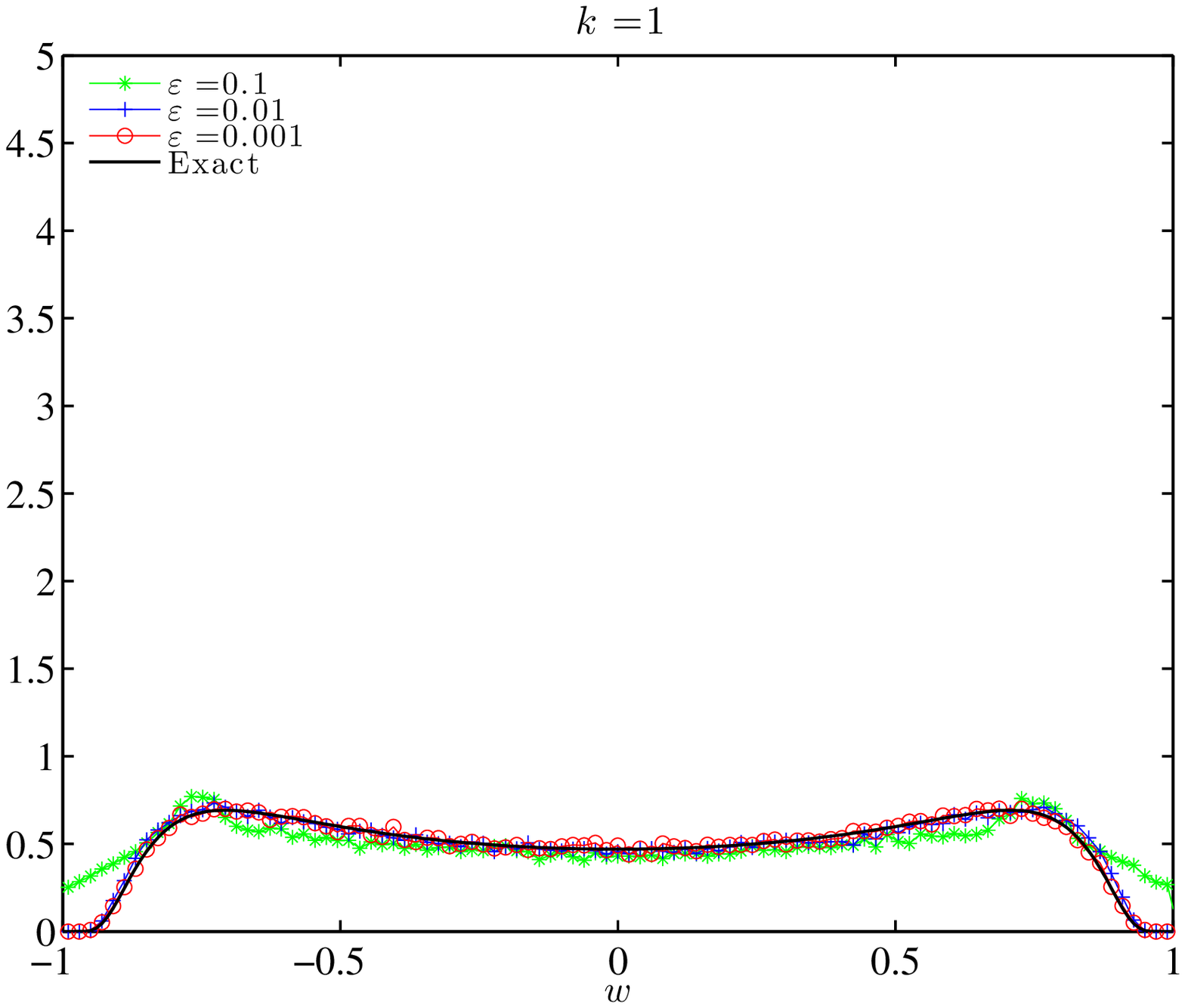}
\\
\includegraphics[scale=0.35]{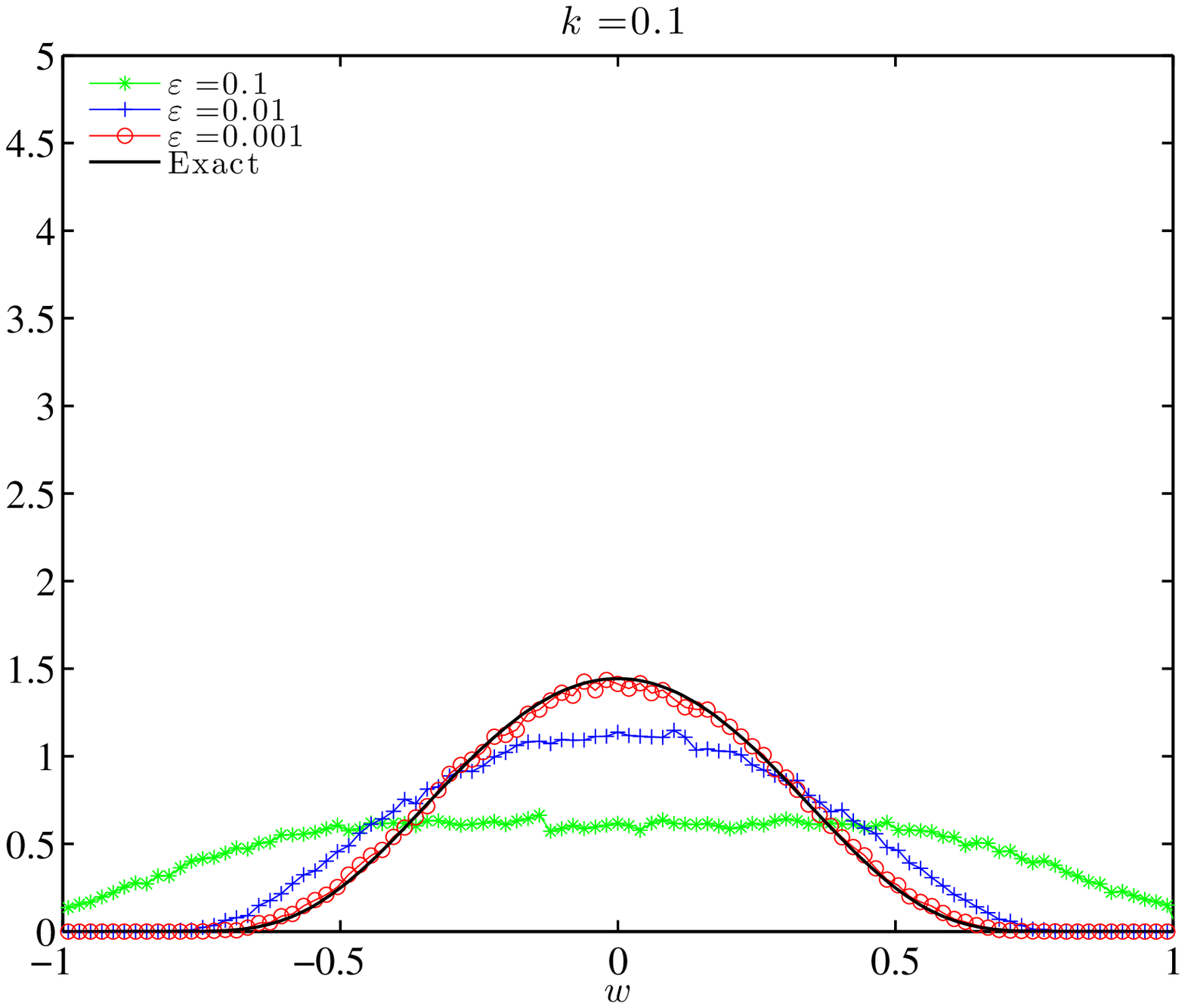}
%&
\includegraphics[scale=0.35]{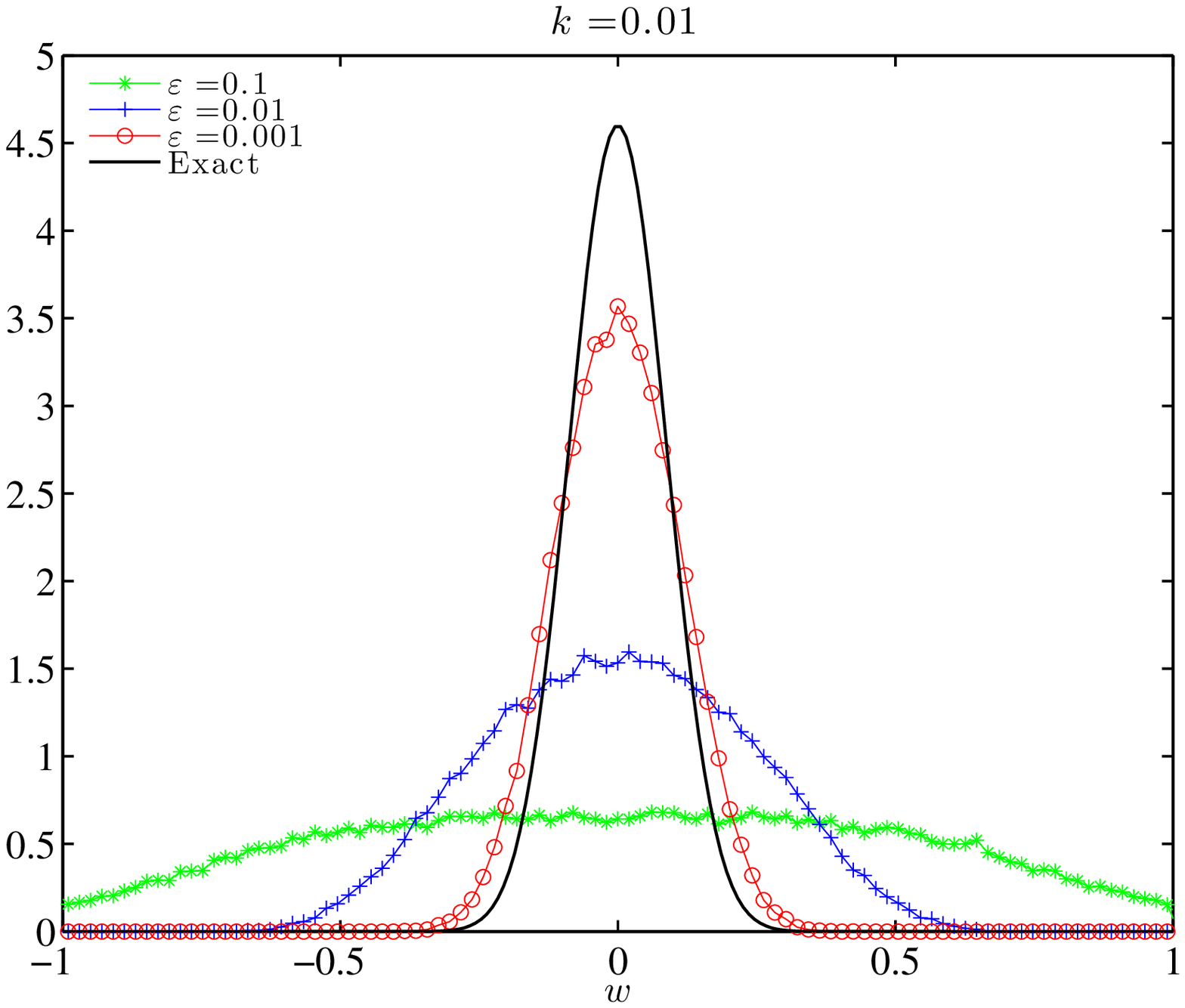}
\caption{Steady solutions of the Boltzmann equation with $P(w,v)=1$ and $D(w)=1-w^2$ in the scaling \eqref{eq:scale} for different values of $\epsi$ and $\varsigma=3$. Continuous lines represent the steady profile of the Fokker--Planck equation. From left to right from top to bottom, we increase the control action, diminishing the value of $\kappa$.}\label{fig:F4}
\end{figure}

In Figure \ref{fig:F4} we simulate the evolution of the probability density $f(w,t)$, using a sample of $N_s=10^5$ agents each of them interacting through the binary dynamic \eqref{eq:BinScale} for different scaling values $\epsi$ and $\Theta$ distributed uniformly on $(-\sigma,\sigma)$, with $\sigma^2=3\epsi\varsigma$, $\varsigma=3$. 
Note that the discrepancy of the steady profiles in Figure \ref{fig:F4} is due to the fact we are simulating the convergence of the Boltzmann equation towards its Fokker-Planck limit. Therefore decreasing $\epsi$ and increasing the size of the sample $N_s$ we can obtain better approximations of the Fokker--Planck profiles.

\begin{figure}[t]
\centering
\includegraphics[scale=0.4]{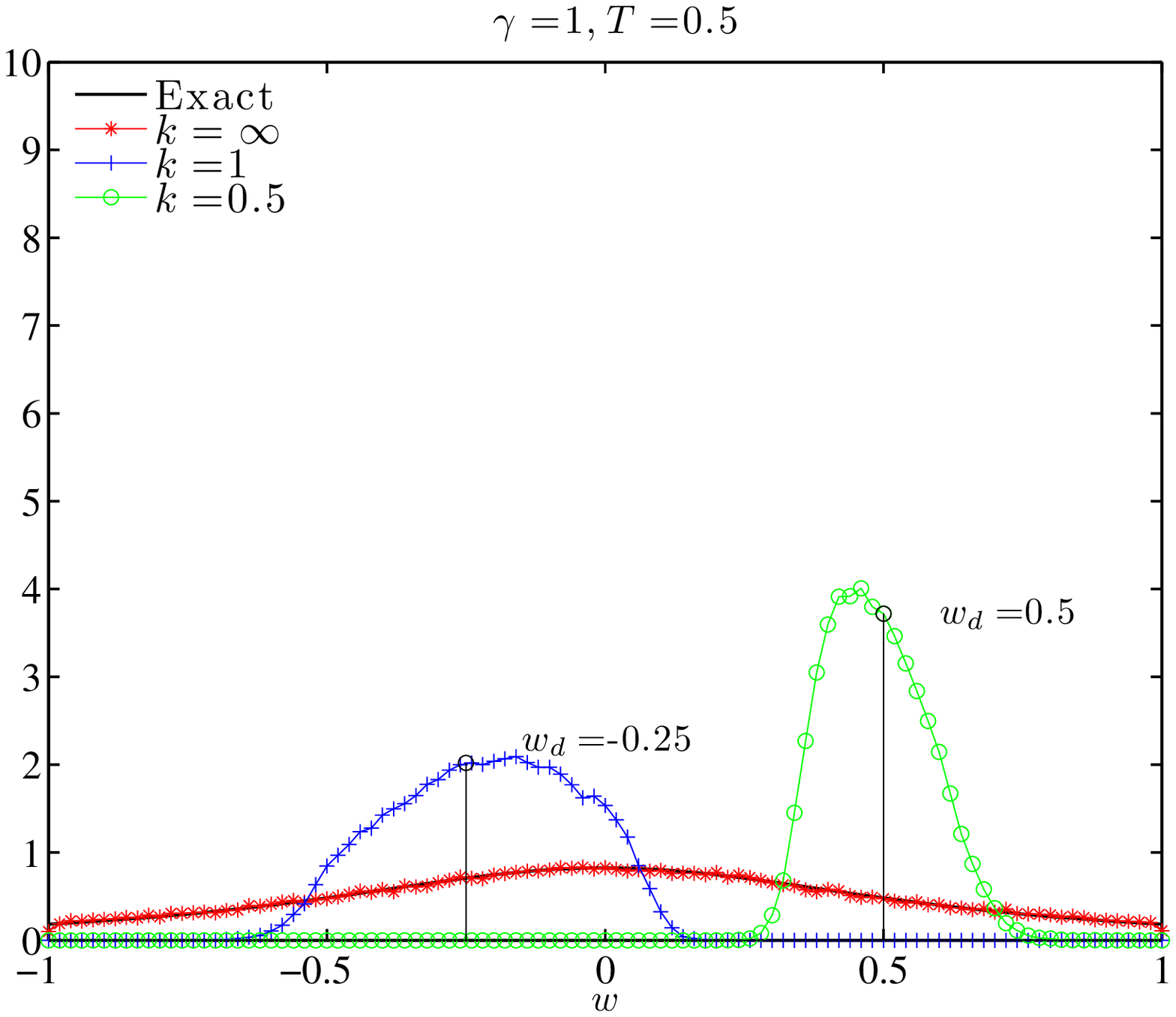}
%&
\includegraphics[scale=0.4]{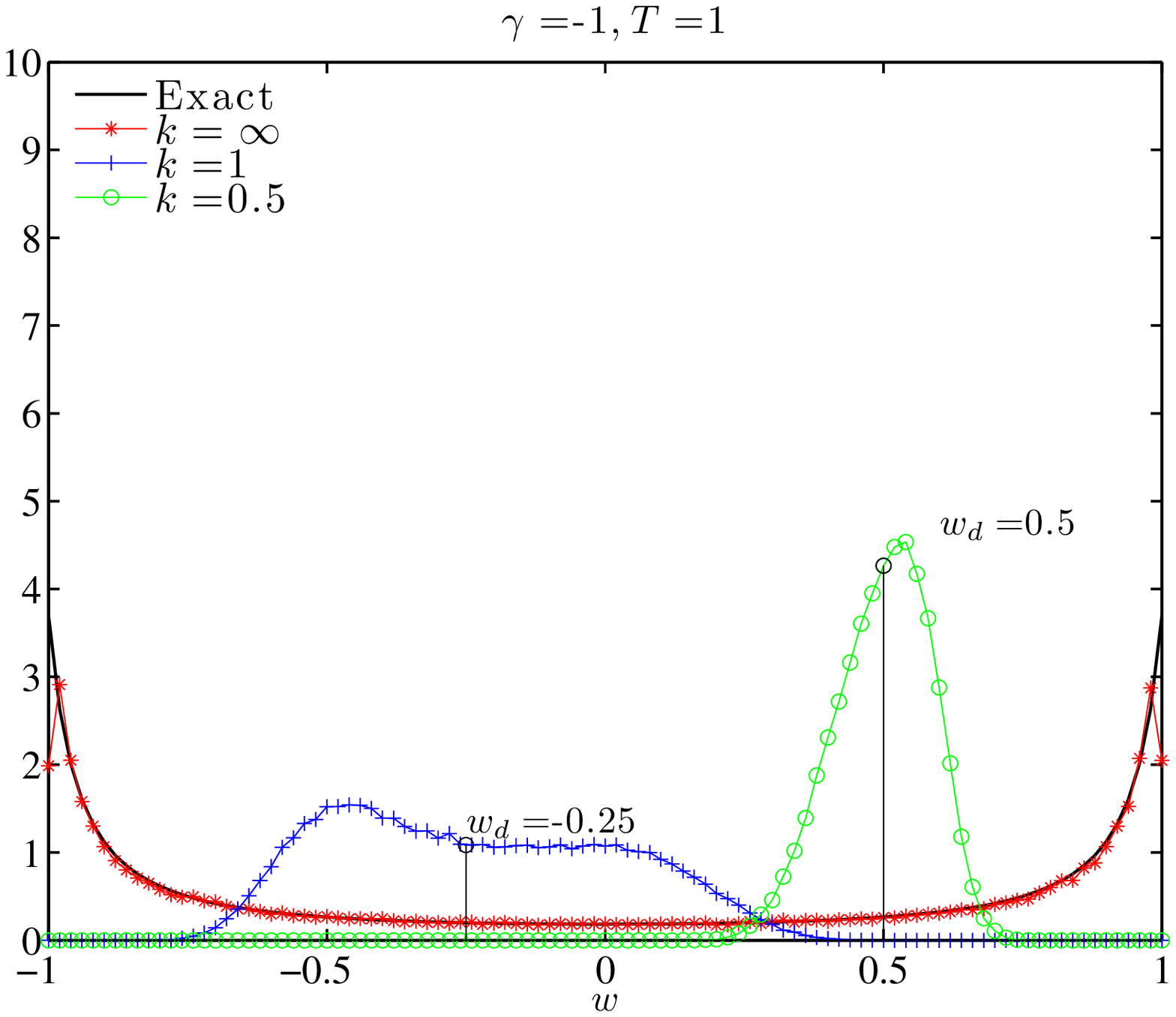}
\\
\includegraphics[scale=0.4]{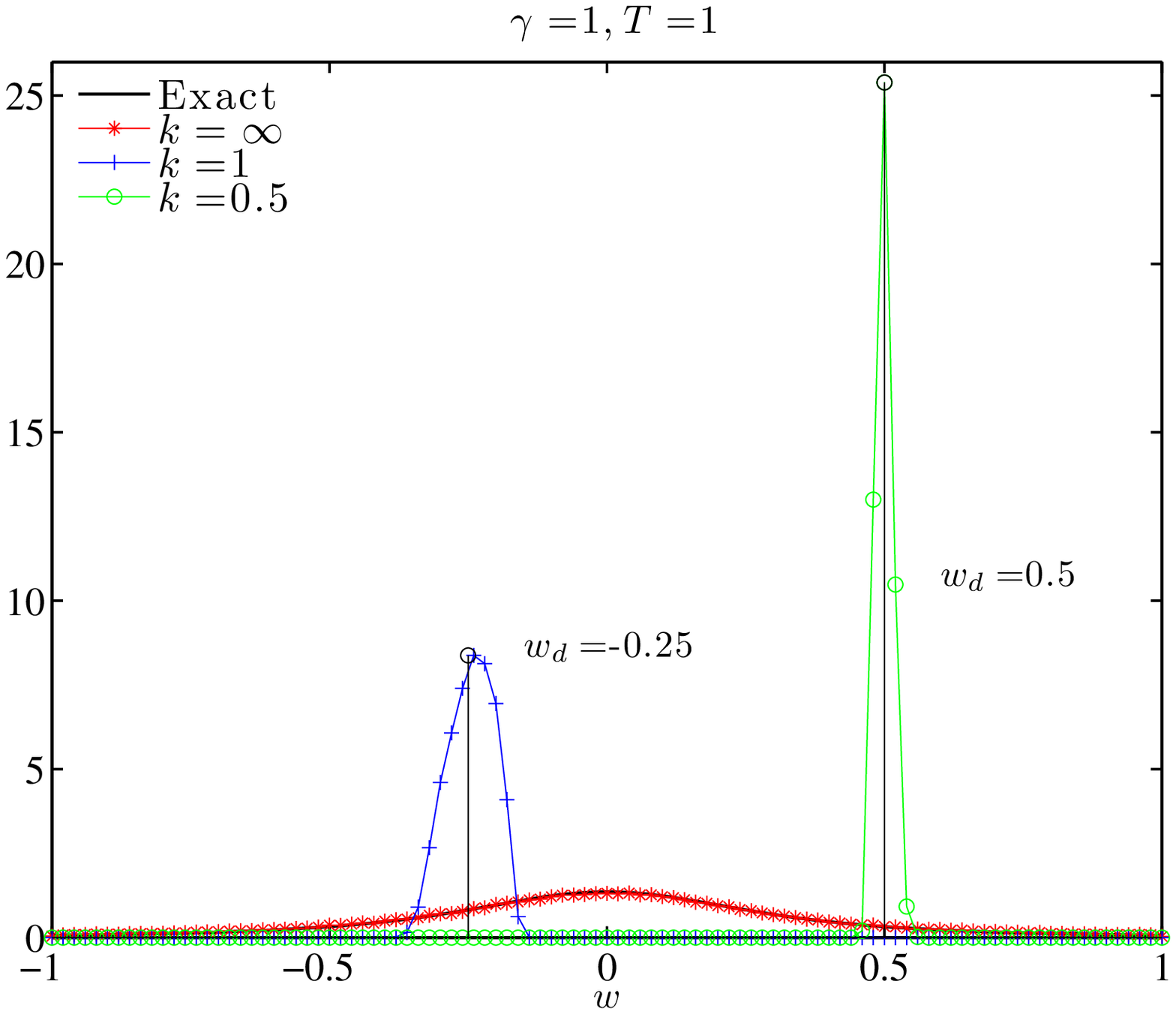}
%&
\includegraphics[scale=0.4]{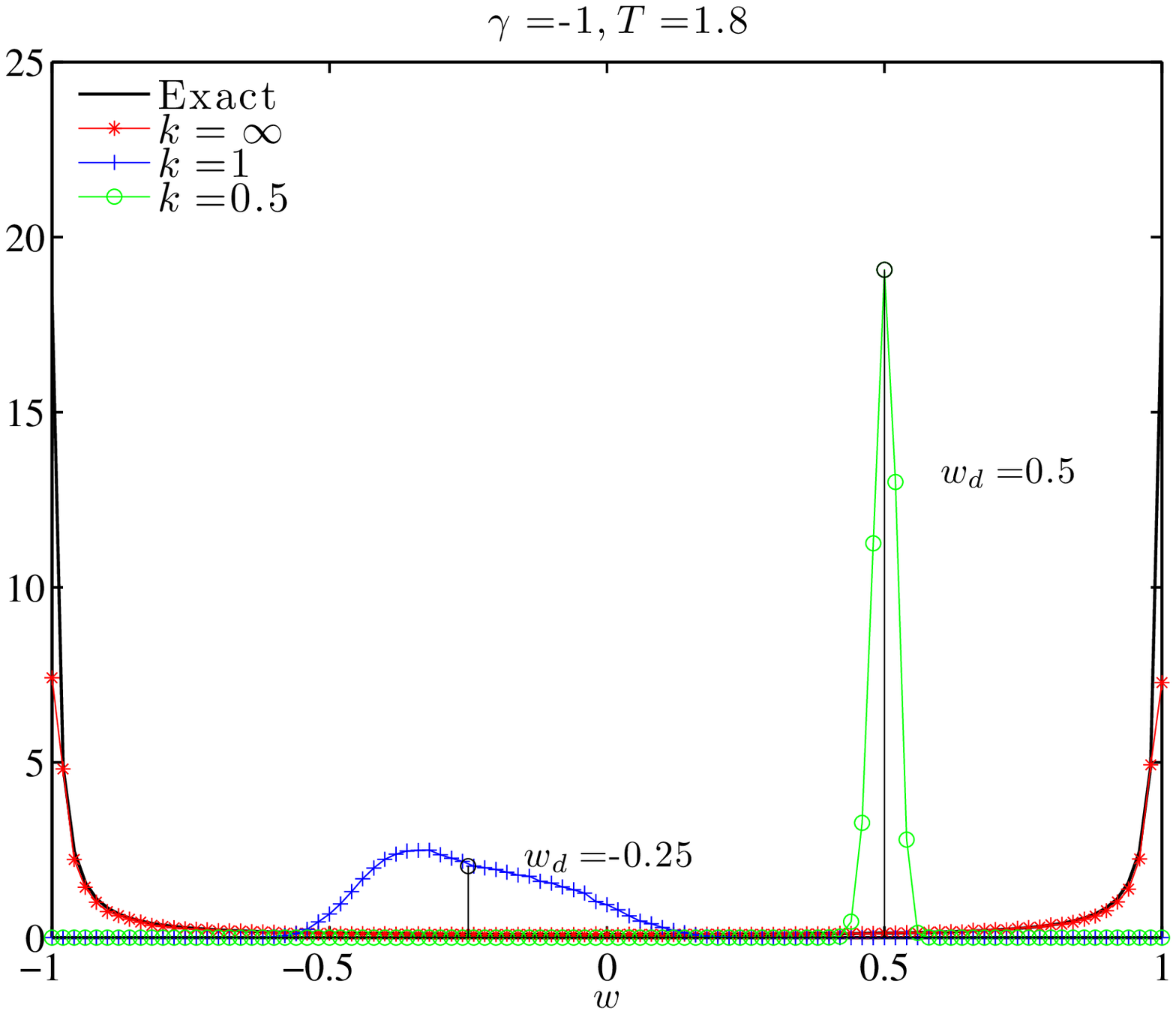}
\caption{Sznajd-type model at different times. The effect of concentration ($\gamma=1$) on the left, and separation ($\gamma=-1$) are visible for the uncontrolled case ($\kappa=\infty$). The action of a mild control $\kappa=1$ and a strong control $\kappa=0.1$ forces the dynamic towards different desired states, respectively  $w_d=-0.25$  and  $w_d=0.5$. As expected the process needs a larger amount of time to control the separation dynamic.}\label{fig:F4b}
\end{figure}

\subsubsection*{Sznajd-type model}
In this test we consider a compromise propensity of the form
\be
P(w,v)=\gamma(1-w^2),\quad \gamma\in\R
\ee 
in absence of diffusion $D(w)=0$. 
Note that, when the initial mean opinion $m(0)=0$, the quasi-invariant opinion limit in absence of control is governed by the mean-field Sznajd's model~\cite{Sznajd00, ANT07}   
\begin{align}\label{eq:MFSz}
\partial_t f =\gamma\partial_w\left(w(1-w^2)f\right).
\end{align}
The model \eqref{eq:MFSz} can be solved explicitly and gives~\cite{ANT07}
\begin{align}\label{eq:eSz}
f(w,t)=\frac{e^{-2\gamma t}}{((1-w^2)e^{-2\gamma t}+w^2)^{3/2}}f_0\left(\frac{w}{((1-w^2)e^{-2\gamma t}+w^2)^{1/2}}\right),
\end{align}
where $f_0(x)$ is the initial distribution.
For $\gamma>0$ we have \emph{concentration} of the profile around zero, conversely for $\gamma<0$ a \emph{separation} phenomena is observed and the distribution tends to concentrate around $w=1$ and $w=-1$.

We simulate the binary dynamic with control corresponding to the above choices starting from an initial mean opinion $m(0)=0$. Our aim is to explore the differences between the controlled concentration and separation dynamics. We choose a scaling parameter $\epsi=0.005$ and a number of sample agents of $N=10^5$. 

In Figure \ref{fig:F4b} we simulate the evolution of $f(w,t)$ for the concentration ($\gamma=1$) and separation ($\gamma=-1$) cases. Starting from the uniform distribution on $\I$, we investigate three different cases: uncontrolled ($\kappa=\infty$), mild control ($\kappa=1$) towards desired state $w_d=-0.25$ and strong control ($\kappa=0.1$) towards $w_d=0.5$. 
The solution profiles in the uncontrolled case, $\kappa=\infty$ coincides with the exact solution profile given by \eqref{eq:eSz}. Observe that separation phenomena implies a slower convergence towards the desired states.

We complete the tests just presented with Table \ref{tab:T4}, where we measure the $L^2$ distance between the average opinion $m$ at final time $T=2$ and the desired state $w_d$, in the separation case, ($\gamma=-1$). We compare the errors for decreasing values of $\kappa$ and for different values of the desired state $w_d$, showing that more effective control implies faster convergence.

\begin{table}[t]
\centering
\begin{tabular}{l c c c c}
\hline\hline
 & $w_d=0.25$ & $w_d=0.5$  & $w_d=0.75$& $w_d=0.95$\\
\hline\hline
 $\kappa=10$  & 1.7139e-01 & 3.428e-01  & 5.1351e-01&  6.5032e-01\\
   \hline
 $\kappa=5$  & 1.1468e-01  & 2.2653e-01& 3.3844e-01 & 4.2362e-01\\
   \hline
 $\kappa=1$  &1.0592e-03 & 1.6027e-03  & 1.5460e-03 & 1.2877e-03\\
   \hline
 $\kappa=0.5$  & 7.0990e-07  & 9.0454e-07  & 6.9543e-07 & 4.9742e-07\\
\hline
\end{tabular}
\caption{$L_2$ distance between $w_d$ and the average opinion $m$ at time $T=2$ for the controlled Sznajd-type model with separation interactions.}\label{tab:T4}
\end{table}

\subsubsection*{Bounded confidence model}
Next, we consider the case of \emph{bounded confidence models}, where the possible interaction between agents depends on the level of confidence they have \cite{HK02, GGL:12}. This can be model through a compromise function which accounts the exchange of opinion only inside a fixed distance $\Delta$ between the agent opinions
\be
P(w,v)=\chi(|w-v|\leq\Delta),
\ee
where $\chi(\,\cdot\,)$ is the indicator function.
\begin{figure}[ht]
\centering
\includegraphics[scale=0.4]{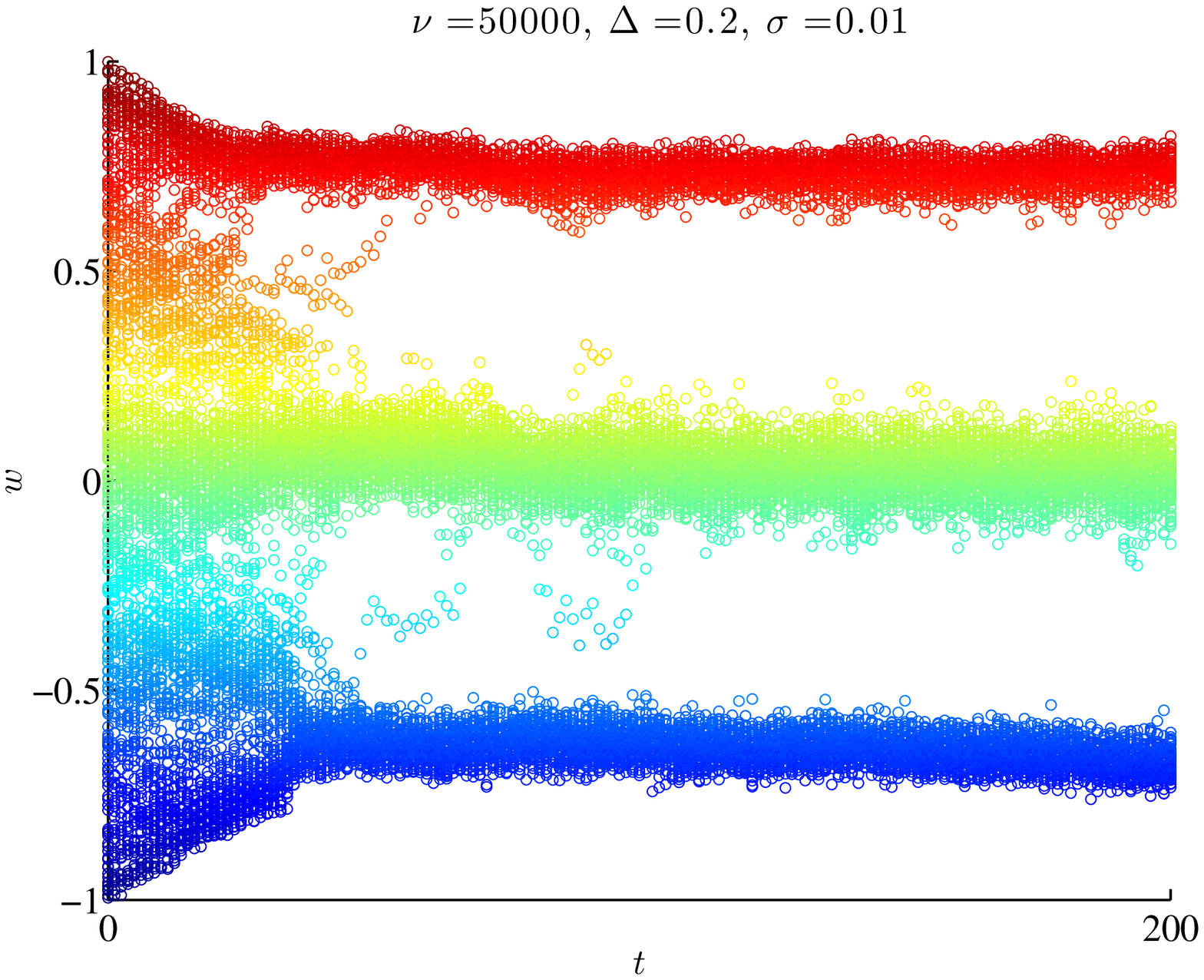}
%&
\includegraphics[scale=0.4]{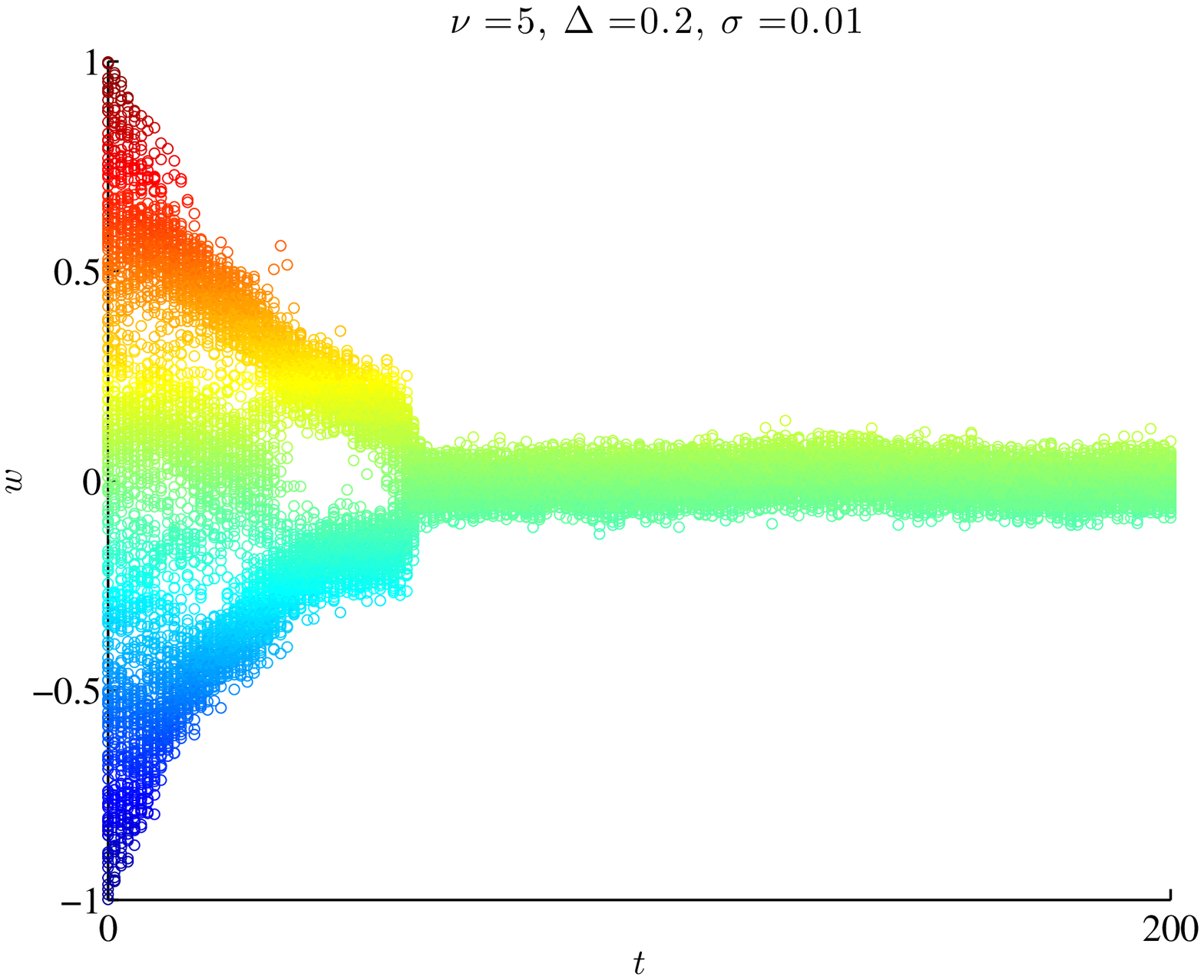}
\\
\includegraphics[scale=0.41]{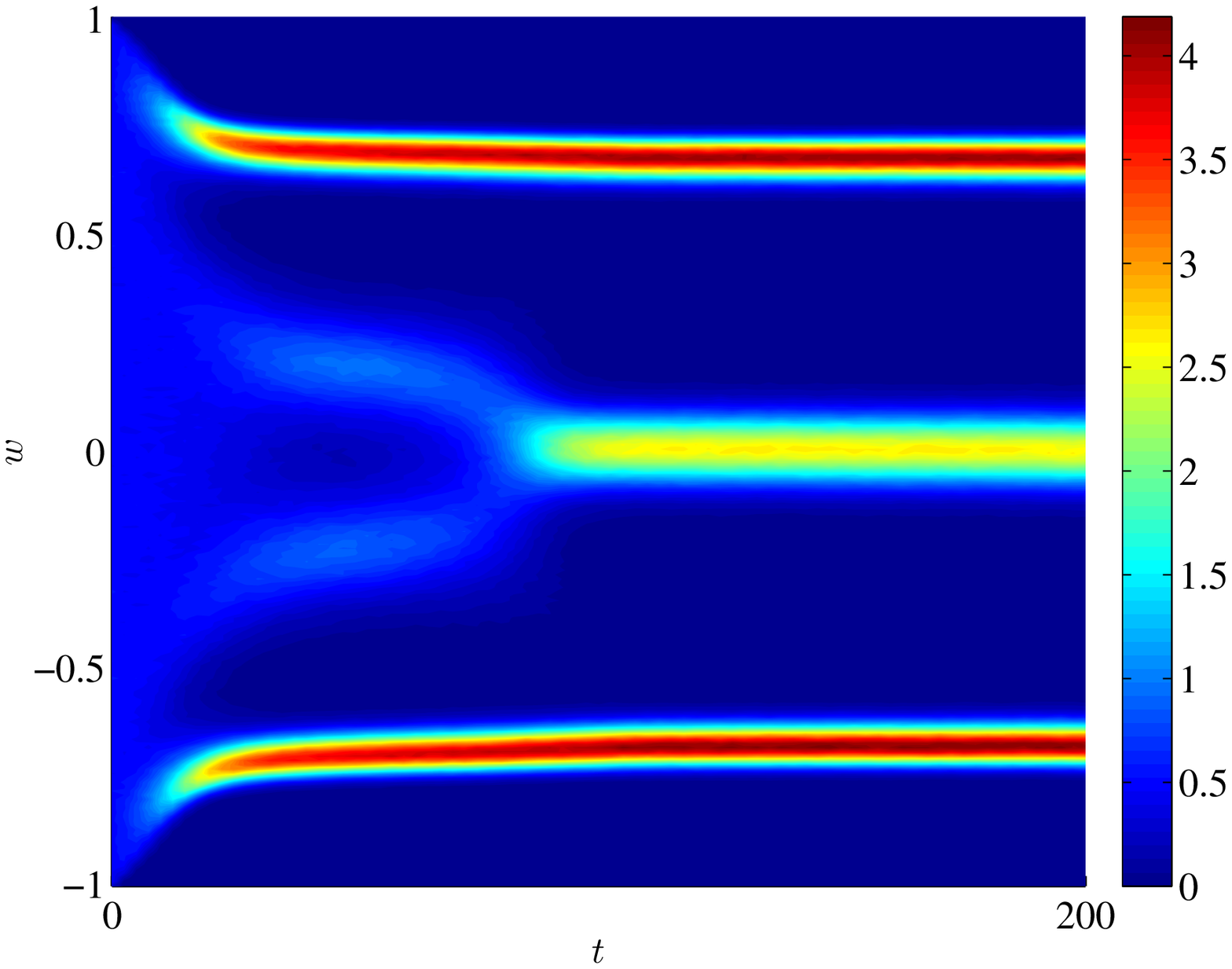}
%&
\includegraphics[scale=0.41]{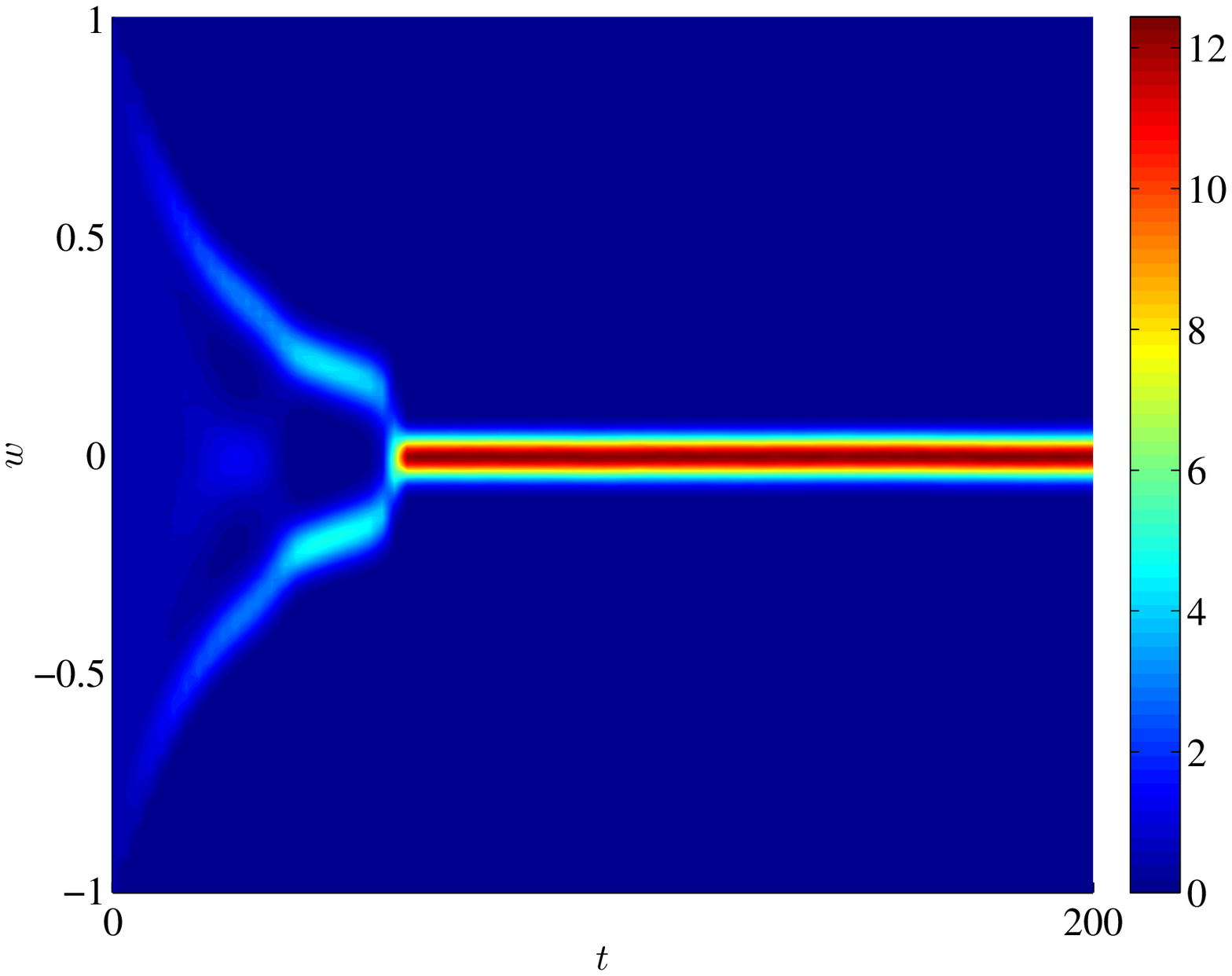}

\caption{Bounded confidence model. On the left the control parameter $\nu=5000$ on the right $\nu=5$. In the top row the result of a particle simulation with $N=200$ agents where the color scale depicts the opinion value. Bottom row represents the evolution of the kinetic density. In both cases the simulation is performed for $\sigma=0.01$ and $\Delta = 0.2$.   }\label{fig:F5}
\end{figure}

In Figure~\ref{fig:F5}, we simulate the dynamic of the agents starting from an uniform distribution of the opinions on the interval $\I=[-1,1]$. The confidence bound is taken $\Delta=0.2$ and the diffusion parameter $\sigma=0.01$. We consider the case without control and with control, letting the system evolve in the time interval $\left[0~ T\right]$, with $T=200$. 
In the left column figures we represents the weak controlled case, with penalization parameter $\nu=5000$, and three mainstream opinions emerge, on the right the presence of the control, $\nu=5$ is able to lead the opinions to concentrate around the desired opinion, $w_d=0$.

Top row of plots shows the evolution of the dynamic at the particle level, with $N=200$. Bottom row represents the same dynamic at the kinetic level, simulation is performed with a sample of $N_s=2\times10^5$ particles with $\epsi=0.05$.

\section{Conclusions}
In this paper we introduced a general way to construct a Boltzmann description of optimal control problems for large systems of interacting agents. The approach has been applied to a constrained microscopic model of opinion formation. The main feature of the method is that, thanks to a model predictive approximation, the control is explicitly embedded in the resulting binary interaction dynamic. In particular in the so-called quasi invariant opinion limit simplified Fokker-Planck models have been derived which admit explicit computations of the steady states. The robustness of the controlled dynamics has been illustrated by several numerical examples which confirm the theoretical results. Different generalizations of the presented approach are possible, like the introduction of the same control dynamic through leaders or the application of this same control methodology to swarming and flocking models. 

\bigskip
\noindent
{\bf Acknowledgements.} This work has been supported
by DFG Cluster of Excellence EXC128, the BMBF KinOpt Project   
 and DAAD 55866082 and by PRIN-MIUR grant {\em ``Advanced numerical methods for
kinetic equations and balance laws with source terms''}.

\bibliographystyle{siam}
\bibliography{BibTex}

\end{document}